\documentclass[12pt]{article}
\usepackage{hyperref,amsfonts,amssymb}
\usepackage{amsfonts,amssymb}
\usepackage{mathtools}

\def\C{\mathbb{C}}
\def\N{\mathbb{N}}
\def\Z{\mathbb{Z}}
\def\Q{\mathbb{Q}}

\def\R{\mathbb{R}}
\def\l{\lambda}

\def\p{ \partial }

\def\bq{ \begin{equation} }
\def\eq{ \end{equation} }
\def\ben{ \begin{eqnarray} }
\def\en{ \end{eqnarray} }

\def\frac#1#2{{#1\over #2}}

\def\on#1#2{\mathop{\vbox{\ialign{##\crcr\noalign{\kern2pt}
$\scriptstyle{#2}$\crcr\noalign{\kern2pt\nointerlineskip}
\kern-2pt$\hfil\displaystyle{#1}\hfil$\crcr}}}\limits}

\begin{document}

\title{Explicit formulas for arithmetic support of differential and difference operators}
\author{Maxim Kontsevich,  Alexander Odesskii}
   \date{}
\vspace{-20mm}
   \maketitle
\vspace{-7mm}
\begin{center}
IHES, 35 route de Chartres, Bures-sur-Yvette, F-91440,
France  \\[1ex]
and \\[1ex]
Brock University, 1812 Sir Isaac Brock Way, St. Catharines, ON, L2S 3A1 Canada\\[1ex]
e-mails: \\%[1ex]
\texttt{maxim@ihes.fr}\\
\texttt{aodesski@brocku.ca}
\end{center}

\medskip

\begin{abstract}

We compute arithmetic support of the formal deformations $D=P+tQ_1+t^2Q_2+...$ of the differential
operator $P=(x\partial_x-r_1)...(x\partial_x-r_k)$, where $r_1,...,r_k\in\Q$ for sufficiently large primes $p$ in terms of the monodromy of $D$ in characteristic zero. An analog of these results is also provided in the case of $q$-difference operators.

\medskip

MSC: 

\medskip

Keywords:  
\end{abstract}

\newpage

\tableofcontents

\newpage

 \section{Introduction}

 \subsection{General description of the problem}

 Consider an algebra of differential operators in one variable $x$ over a field of characteristic zero, for example, $\Q[x,\partial_x]$. It does not have non-trivial 
 finite-dimensional representations\footnote{Indeed, let $x\mapsto U,~\partial_x\mapsto V$ be an $n$-dimensional matrix representation, $n\ge 1$. Then $VU-UV=1$. By taking the trace of both sides, we obtain $0=n$, which is a contradiction in characteristic zero. However, it is not a contradiction in positive characteristic.}. On the other hand, an algebra of differential operators over 
 a field of characteristic $p$, such as $\Z/(p)[x,\partial_x]$ possesses a two-parametric family of $p$-dimensional matrix representations $\rho_{p,\xi,\eta}$ given by 
 \begin{equation}\label{rhop}\rho_{p,\xi,\eta}(x)=
		\begin{pmatrix}
			\xi & 0 & & 0 \\
			1 & \xi & \ddots & \\
			& \ddots & \ddots & 0 \\
			0 & & 1 & \xi
		\end{pmatrix}, \quad
        \rho_{p,\xi,\eta}(\partial_x)=
		\begin{pmatrix}
			\eta & 1 & & 0 \\
			0 & \eta & \ddots & \\
			& \ddots & \ddots & p-1 \\
			0 & & 0 & \eta
		\end{pmatrix}.
	\end{equation}
 % \begin{equation}\label{rhop}\rho_{p,\xi,\eta}(x)=\left( \begin{array}{cccccccccc}
%\xi &  0 & & 0\\
% 1 &  \xi &  0 &  &\\
%&  1 &  \xi &   0 \\
% &  & \ddots &  \ddots  &   \\
 
% 0&&~~~  1& \xi\end{array}  \right) , ~~~ \rho_{p,\xi,\eta}%(\partial_x)=\left( \begin{array}{cccccccccc}
%\eta &  1 & & 0\\
 %0 &  \eta &  2 &  &\\
%&  0 &  \eta &   3 \\
% &  & \ddots &  \ddots  &   \\
 
% 0&&  0& \eta\end{array}  \right).
% \end{equation}

For a given differential operator $P\in\Z[x,\partial_x]$ and all primes $p=2,3,5,...$ we can reduce $P$ modulo $p$ and compute 
determinant\footnote{One can show that in fact $\det\rho_{p,\xi,\eta}(P)\in\Z/(p)[\xi^p,\eta^p].$} 
$$\det\rho_{p,\xi,\eta}(P)\in\Z/(p)[\xi,\eta].$$
It would be interesting to understand how these polynomials in $\xi,\eta$ with
coefficients in $\Z/(p)$ depend on $p$. It turns out that for a generic
differential operator $P\in\Z[x,\partial_x]$ the polynomials $\det\rho_{p,\xi,\eta}(P)$ depend on $p$ rather chaotically and do not possess any explicit
formula. On the other hand, if $P$ is a so-called motivic operator, then for all but finitely many $p$ there
exists a simple universal formula for $\det\rho_{p,\xi,\eta}(P)$, and the existence
of such a formula can be taken as an alternative definition of motivic
differential operators with rational coefficients, see \cite{Bostan}, \cite{K1}. In this paper, we study the
determinants of formal deformations of a simplest class of motivic operators in representations
$\rho_{p,\xi,\eta}$ and find explicit formulas for them in terms of monodromy
of formal solutions in characteristic zero. 

A similar problem can be posed in the case of an algebra of $q$-difference operators $\Q\langle q,x,y\rangle$ with generators $q,x,y$ and commutation relations $yx=qxy$, where
$q$ is a central element.  If $q$ is a primitive $n$-th root of unity, then
this algebra has a family of $n$-dimensional representations $\rho_{n,\xi,\eta}$ given by (\ref{rhon}). Therefore, for a $q$-difference operator $P\in \Q\langle q,x,y\rangle$ we
obtain a series of polynomials 
$$\det \rho_{n,\xi,\eta}(P)\in \Q[q,\xi,\eta]/\Phi_n(q),~~~n=1,2,3,...$$
where $\Phi_n(q)\in\Z[q]$ stands for $n$-th cyclotomic polynomial. There
exists a class of $q$-difference operators which can be called motivic by similarity with the case of differential operators, and such that these
determinants are given by very simple formulas. In this paper, we study the
determinants of formal deformations of motivic $q$-difference operators of the
form 
$$P=(y-q^{l_1})...(y-q^{l_k}),~~~~l_1,...,l_k\in\Z.$$

\subsection{Motivation and more detailed description}

 Let us explain the original motivation for the study of determinants of 
 differential (resp. $q$-difference) operators in representations $\rho_{p,\xi,\eta}$, $p=2,3,5,...$ (resp. its $q$-analogs $\rho_{n,\xi,\eta}$, $n=1,2,3,...$). Notice that in our previous paper \cite{KO} 
we studied these determinants for differential operators in the particular
case $\xi=\eta=0$ using a different and more technical approach, which does
not work in the general case. 

The idea of a motive was introduced by A. Grothendieck in an attempt to define a universal cohomology theory for algebraic varieties. For families of algebraic varieties, one should get a kind of motivic sheaf on the base of
the family. In particular, in the de Rham realization one obtains an algebraic
vector bundle with flat connection (the so-called Gauss-Manin connection) or, 
in classical terminology, Picard-Fuchs differential equations. Solutions of
such equations can be written as integrals of algebraic functions in several
variables depending on the parameters. There are several hypothetical characterizations of motivic connections (also called connections of geometric origin). One of such characterizations is the so-called \textit{ global nilpotence}, which is the vanishing of so-called $p$-curvature for sufficiently large primes $p$.

Let $X/R$ be a smooth algebraic variety over $\text{Spec}\ R$ where $R$ is a finitely generated ring without zero divisors and assume that $\Z\subset R$, i.e. the field of fractions of $R$ has characteristic zero. Let $\mathcal{E}$ be a vector bundle
on $X$ with flat connection $\nabla$. Then for all primes large enough $p$ the variety
$X$ and the space $\text{spec}R$ have good reduction modulo $p$, and we have a well defined $p$-curvature of connection $\nabla$ mod $p$. Recall that the $p$-curvature operator corresponding to a vector field $\mu\in\Gamma(X,T_X)$ is 
defined by $(\nabla_{\mu})^p-\nabla_{\mu^p}$ mod $p$. 

N. Katz proved \cite{katz} that the flat connection of algebro-geometric origin, or in other words motivic, that is, the Gauss-Manin connection on the de Rham cohomology of a family of
varieties, and its subquotients, is globally nilpotent, that is, as we said, has a nilpotent $p$-curvature for all sufficiently large $p$. Global nilpotence is conjectured to be equivalent to several other arithmetic properties, including the property of solution to be a $G$-series in the sense of Siegel (see \cite{A}, \cite{AB}), and the latter are expected to be of geometric origin by Bombieri-Dwork conjecture.

%As a partial evidence, he proved that the nilpotence implies two properties of 
%connections of geometric origin

%1) they have regular singularities

%2) the monodromy around any divisor is unipotent.

\

In \cite{K1} the notion of arithmetic support of a holonomic $\mathcal{D}$-module was introduced and a generalization of geometricity conjecture to $\mathcal{D}$-modules with irregular singularities was suggested\footnote{An important property of the arithmetic support that it is Lagrangian was proven 
in \cite{Bit} in the case of differential operators, and in \cite{Et} in the case of $q$-difference operators.}. Moreover, one can reduce this generalization to a question concerning polynomial differential operators in one variable.

 Let 
 \begin{equation}\label{P0}
 P=\sum_{0\leqslant i+j\leqslant N} a_{i,j}x^i\partial_x^j
 \end{equation}
 where $a_{i,j}\in R$. Then for any prime $p=2,3,5,...$ we have a determinant 
 $$\det \rho_{p,\xi,\eta}(P)=\sum_{0\leqslant i+j\leqslant N}b_{i,j}\xi^{pi}\eta^{pj}~~\in~ R/(p)[\xi^p,\eta^p],$$
 see \cite{K1} for details. Here, $\rho_{p,\xi,\eta}$ is a $p$-dimensional representation defined by (\ref{rhop}). The above polynomial is in a sense the characteristic polynomial of $p$-curvature.

 The analog of the global nilponence condition is that for all $p\gg 2$ we have 
 \begin{equation}\label{bc}
 b_{i,j}=c_{i,j}^p~~~ \text{mod} ~p
 \end{equation}
 for some $c_{i,j}\in R\otimes\Q$ independent of 
 $p$. The following conjecture was suggested in \cite{K1}: 
 
 {\it $P$ satisfies the above property iff $T$ is an extension of subquotients of exponentially-motivic $\mathcal{D}$-modules on $\mathbb{A}^1$.}

 In particular, if a motivic operator $P$ given by (\ref{P0}) has rational 
 coefficients $a_{i,j}\in\Q$, then 
 \begin{equation}\label{fp}
 \det \rho_{p,\xi,\eta}(P)=\sum_{0\leqslant i+j\leqslant N}c_{i,j}\xi^{pi}\eta^{pj}~~\text{mod}~p
 \end{equation}
 for $p\gg 2$, where $c_{i,j}\in\Q$ are independent of $p$.

 Exponentially motivic $\mathcal{D}$-modules are such that all solutions of the corresponding system of differential equations can be written as integrals of products of algebraic function and exponential of algebraic function, see \cite{K1}.

 In the case of $\mathcal{D}$-modules with regular singularities, one gets just a reformulation of the global nilpotence criterion of geometricity. 
 In the article \cite{Bostan} several examples of algebraic linear equations originating from statistical physics were discussed, in which such a criterion can be applied. In addition, in the same article an example of a $\mathcal{D}$-module with irregular singularity was considered, where the characteristic polynomial of $p$-curvature lifts to a polynomial with coefficients in $\mathbb Q$, i.e., satisfies the same generalization of the global nilpotence criterion as in \cite{K1}.

 A natural question arises: consider a formal deformation 
 $$D=P+tQ_1+t^2Q_2+...$$
 of an operator $P$ satisfying the condition (\ref{bc}). What can we say about the expansion coefficients of $\det \rho_{p,\xi,\eta}(D)$ with respect to the deformation parameter $t$?

 In this paper, we answer this question in the simplest case 
 $$P=\prod_{i=1}^k(x\partial_x-r_i),~~~r_i\in\Q.$$

 We also consider a $q$-analog of the determinants described above, where differential operators are replaced by $q$-difference operators, and prime numbers are replaced by cyclotomic polynomials $\Phi_n(q)\in\Z[q]$, $n=1,2,...$, and $p$-dimensional representations $\rho_{p,\xi,\eta}$ are replaced by $n$-dimensional representations $\rho_{n,\xi,\eta}$ given by (\ref{rhon}). Notice that a $q$-difference version of the generalized geometricity conjecture was also proposed in \cite{K1}. In this paper, we study the determinants of the deformations of a family\footnote{Notice that a $q$-analog of our results can also be done in the case $l_i\in\Q$ but we restrict ourselves to the integer case in order to simplify the notation.} of difference operators
 $$P=\prod_{i=1}^k(y-q^{l_i}),~~~l_i\in \Z$$
 where $y$ is a $q$-shift operator $yf(x)=f(qx)$. 

 Let $P$ be a differential (resp. difference) operator with a simple explicit formula for $\det\rho_{p,\xi,\eta}(P)$ (resp. $\det\rho_{n,\xi,\eta}(P)$). In this 
 case we can write our determinants in terms of traces by
 $$\det(P+Q)=\det(P)\cdot \det(1+P^{-1}Q)=\det(P)\cdot\exp(\text{tr}\log(1+P^{-1}Q))$$
 where we omit our notation for $p$ or $n$ dimensional representations, 
 to simplify notations. This means that it is natural to extend our algebra 
 of differential or difference operators by adding two-sided inverse $P^{-1}$
 and study quotient by commutant for this extended algebra. For example, 
 let $P=\prod_i(x\partial_x-r_i),~r_i\in\Z$, then the quotient of the extended 
 algebra by commutant has a basis $\{(x\partial_x)^{-i};~i=1,2,...\}$. In the case of difference operators, the corresponding commutant has a basis $\{1,(y-1)^{-1},(y-1)^{-2},...\}$.

 \subsection{Content of the paper}

 In Section 2 we find an explicit formula for traces of $p\times p$ matrices 
 $$\rho_{p,\xi,\eta}(P^{-1}D_1P^{-1}D_2...D_NP^{-1}),~~~p\gg 2$$
 where $D_1,...,D_N\in R[x^{\pm 1},\partial_x]$ are differential operators over
 a finitely generated subring of a field characteristic zero, and $P$ is a 
 monic polynomial in $x\partial_x$ with rational roots. Based on this result, we compute the determinants of formal deformations 
 $$P+tQ_1+t^2Q_2+...$$
 in representations $\rho_{p,\xi,\eta}$. See formula (\ref{trA}) for the traces and formula (\ref{detD}) for the determinants. Notice that in the r.h.s. of both formulas we have expressions in characteristic zero independent of $p$ which is equal to the l.h.s. modulo $p$.

 In Section 3 we investigate formal solutions of the equation
 $$(P+tQ_1+t^2Q_2+...)f=0$$
 which are deformations of $f_0$, a solution of $Pf_0=0$. Here $P$ is a monic polynomial in $x\partial_x$ of degree $k$. We prove that the space of such solutions is $k$-dimensional over formal series in $t$. We also define (after an extension of scalars) a basis in the space of solutions such that its elements can be seen as algebraic analogs of generalized eigenvectors of the monodromy operator. We also write 
 our determinants in terms of algebraic analog of eigenvalues of the monodromy operator.

 In Section 4 we explain the relation with the monodromy operator in an analytic situation.

 In Section 5 we give a $q$-analog of previous results.

 In Section 6 we show that in the simplest nonzero case $P=1$, $k=0$ there exist more precise formulas for certain classes of deformations of $P=1$ by $q$-difference operators. 

 In Section 7 we present some results based on computer  experiments about adic properties of coefficients of $t$-expansions of our determinants 
 for certain examples.

 \subsection{Notations}

 Let $K$ be an arbitrary ring\footnote{All rings in this paper are associative and unital.}. In the paper we use standard notations for various rings over $K$, such as $K[u]$ for polynomials, 
 $K[[u]]$ for formal power series, $K(\!(u)\!)$ for formal Laurent series, $K(u)$ 
 for rational functions in $u$, and $u$ might be an expression in other 
 variables.

We will often iterate the above constructions of rings, and the result may depend on 
the order of iterations. For example, we have
$$K[u_1][[u_2]]\not\cong K[[u_2]][u_1]$$
which will often be important in our constructions and results.

We also use the notation $K[x,\partial_x]$ for the ring of differential operators in $x$ with coefficients in $K$
$$K[x,\partial_x]=\Big\{\sum_{i,j=0}^ma_{i,j}x^i\partial_x^j,~~~a_{i,j}\in K\Big\}.$$
Moreover, if $v_1,...,v_m\in K[x,\partial_x]$, then we denote\footnote{We hope that this abuse of notations will not create any misunderstanding.} by
$$K[x,\partial_x,v_1^{-1},...,v_m^{-1}]$$
the noncommutive ring over $K$ generated by $x,\partial_x,v_1^{-1},...,v_m^{-1}$ with defining relations
$$\partial_xx-x\partial_x=1,~v_1v_1^{-1}=v_1^{-1}v_1=...=v_mv_m^{-1}=v_m^{-1}v_m=1.$$
In other words, we extend $K[x,\partial_x]$ by adding two-side inverses of 
the elements $v_1,...,v_m\in K[x,\partial_x]$. In particular,
$$K[x^{\pm 1},\partial_x]=\Big\{\sum_{i=-m}^m\sum_{j=0}^na_{i,j}x^i\partial_x^j,~~~a_{i,j}\in K\Big\}.$$

 \section{Explicit formulas for traces and determinants of differential operators}

  Let ${\bf k}$ be a field of characteristic zero. Let ${\bf k}[x^{\pm 1},\partial_x]$ be the ring of differential operators over ${\bf k}$ on the punctured line $\mathbb{A}^1_{\bf k}\setminus\{0\}$. We extend ${\bf k}[x^{\pm 1},\partial_x]$ by inverting arbitrary monic polynomials in ${\bf k}[x\partial_x]$ and denote this extension by $\overline{{\bf k}[x^{\pm 1},\partial_x]}$. Notice that any element $S\in\overline{{\bf k}[x^{\pm 1},\partial_x]}$ can be written uniquely in the form \begin{equation}\label{S}
  S=\sum_{j=l}^mx^jf_j(x\partial_x)
  \end{equation}
  where $l,m\in\Z$ and $f_j(x\partial_x)\in {\bf k}(x\partial_x)$ are rational functions. For multiplication of such expressions one can use commutation  relations:
  \begin{equation} \label{comrel}
  f(x\partial_x)x^j=x^jf(x\partial_x+j),~~~j\in\Z
  \end{equation}
  for an arbitrary rational function $f(x\partial_x)$.

  Let $R\subset {\bf k}$ be a subring. Assume that $R$ is finitely generated over $\Z$. Notice that for any nonzero element $s\in R$ its reductions $s~\text{mod}~p\in R/(p)$ are non-zero for all primes large enough $p$.  Fix
  rational numbers $r_1,\dots,r_k\in\Q$ and assume that these numbers 
  belong\footnote{We can always achieve this by extending the ring $R$.} to 
  $R$.  Let $\text{Diff}=R[x^{\pm 1},\partial_x]$ be the ring of differential operators over $R$, and $\overline{\text{Diff}}\subset \overline{{\bf k}[x^{\pm 1},\partial_x]}$ be a subring consisting of expressions (\ref{S}) such that the numerators of all $f_j(x\partial_x)$ belong to $\text{Diff}$ and the denominators are products of linear polynomials of the form $x\partial_x-r_i+m,~~~1\leqslant i\leqslant k,~m\in\Z$. In other words, 
  $$\overline{\text{Diff}}=R[x^{\pm 1},\partial_x,(x\partial_x-r_i+m)^{-1},~i=1,...,k,~m\in\Z].$$
  Notice that $\text{Diff}\subset \overline{\text{Diff}}$.

  For each large enough prime number $p$ and formal parameters $\xi,~\eta$ we define a $p-$dimensional representation $\rho_{p,\xi,\eta}$ of the ring $\overline{\text{Diff}}$ in the vector space $R/(p)(\xi,\eta)[x]/(x^p)$ by 
  $$\rho_{p,\xi,\eta}(x)=x+\xi,~~~\rho_{p,\xi,\eta}(\partial_x)=\partial_x+\eta.$$
  In the basis $\{1,x,...,x^{p-1}\}$ these operators act by $p\times p$ matrices with coefficients in $R/(p)(\xi,\eta)$ given by (\ref{rhop}).
  
  For a formal parameter $\varepsilon$ we define an infinite-dimensional representation $\rho_{\varepsilon}$ of the ring $\overline{\text{Diff}}$ in the vector space over ${\bf k}(\varepsilon)$ with a basis $\{x^{\varepsilon+m},~m\in\Z\}$ by
  
  \begin{equation}
  \begin{aligned} \label{repep}
    \rho_{\varepsilon}(x):~x^{\varepsilon+m}\mapsto x^{\varepsilon+m+1},~~~~~~~ \\
    \rho_{\varepsilon}(\partial_x):~x^{\varepsilon+m}\mapsto (\varepsilon+m)x^{\varepsilon+m-1}.
  \end{aligned}
\end{equation}

  Let ${\bf \bar{k}}$ be an algebraic closure of the field ${\bf k}$. Define a linear operator $T_{\varepsilon}:~{\bf \bar{k}}(\varepsilon)\to {\bf \bar{k}}(\varepsilon)$ by \begin{equation}\label{T}
  T_{\varepsilon}(\varepsilon^i)=0,~i\geqslant 0,~~~~~~~T_{\varepsilon}\Bigg(\frac{1}{(\varepsilon-c)^j}\Bigg)=\frac{1}{\varepsilon^j},~j>0,~c\in {\bf \bar{k}}
  \end{equation}
  Notice that $T_{\varepsilon}$ preserves ${\bf k}(\varepsilon)$ because it is invariant with respect to the Galois group $Gal({\bf \bar{k}}/{\bf k})$.

  {\bf Theorem 2.1.} Let $A=P^{-1}D_1P^{-1}D_2...D_NP^{-1}\in\overline{\text{Diff}}$ where \begin{equation}\label{P}P=\prod_{j=1}^k(x\partial_x-r_j),~~~r_j\in\Q
  \end{equation}
   and $D_1,...,D_N\in \text{Diff}$. Then for all large enough prime numbers $p>C$ where $C$ depends only on $A$  we have
  \begin{equation}\label{trA}\text{tr}~ (\rho_{p,\xi,\eta}(A))=T_{\varepsilon} \text{Coef}_{x^{\varepsilon}}(\rho_{\varepsilon}(A)x^{\varepsilon}) |_{\varepsilon=-\xi^p\eta^p}~~~\text{mod}~p.
  \end{equation}

  {\bf Proof.} 
  Each $D_j$ can be written in the form $D_j=\sum_{i=-d_j}^{d_j}x^if_{i,j}(x\partial_x)$ where $f_{i,j}$ are polynomials. Using commutation relations (\ref{comrel}) we can write our operator $A$ in the form \begin{equation}\label{A}
  A=\sum_{j=-d_1-...-d_N}^{d_1+...+d_N}x^jg_j(x\partial_x)
  \end{equation}
  where $g_j(x\partial_x)\in\overline{\text{Diff}}$ are rational functions with coefficients in $\Q\subset {\bf k}$. We need the following

  {\bf Lemma 2.1.} The following identities hold
  \begin{equation}\label{trf1}
  \text{tr}~(\rho_{p,\xi,\eta}(x^i(x\partial_x)^j))=0~~~\text{if}~~~i\in\Z~~~\text{and}~~~0\leqslant j\leqslant p-2,
  \end{equation}
   \begin{equation} \label{tr}
  \text{tr}~(\rho_{p,\xi, \eta}(x^i(x\partial_x-r)^{-j}))=\begin{cases}0, &~ i \ne 0~~~\text{mod}~~~p \\  \frac{1}{(-\xi^p\eta^p)^j}, &~~~i=0~\text{and}~0<j\leqslant p \end{cases}
  \end{equation}
  where $r\in\Q$.   
  
  {\bf Proof.} Let $\lambda$ be a formal parameter. Let us prove the following more general identities:
  \begin{equation}\label{detf}
  \det(\rho_{p,\xi,\eta}(x\partial_x-\lambda))=\xi^p\eta^p+\lambda-\lambda^p
  \end{equation}

  \begin{equation}\label{trf2}
  Tr(\rho_{p,\xi, \eta}(x^i(x\partial_x-\lambda)^{-j}))=\begin{cases}0, &~ i \ne 0~~~\text{mod}~~~p \\  \frac{1}{(-\xi^p\eta^p-\lambda+\lambda^p)^j}, &~~~i=0~\text{and}~0<j\leqslant p \end{cases}
  \end{equation}
  If we set $\lambda=r\in\Q$, then $r^p=r$ mod $p$ if $p$ is larger than the denominator of $r$ and we get the identity (\ref{tr}).

   Let $\phi(\lambda)=-\det(\rho_{p,\xi,\eta}(x\partial_x-\lambda))\in\Z/(p)[\xi,\eta,\lambda]$. From the definitions it follows that $\phi(\lambda)$ is a monic polynomial in $\lambda$ of degree $p$. We have 
  $$\phi(0)=-\det(\rho_{p,\xi,\eta}(x\partial_x))=-\det(\rho_{p,\xi,\eta}(x))\det(\rho_{p,\xi,\eta}(\partial_x))=-\xi^p\eta^p$$
  because $\rho_{p,\xi,\eta}(x)$ and $\rho_{p,\xi,\eta}(\partial_x)$ are triangular $p\times p$ matrices with $\xi$ and $\eta$ on diagonals. We also have $x\cdot (x\partial_x-\lambda)\cdot x^{-1}=x\partial_x-1-\lambda$, see (\ref{comrel}). Taking a trace of this relation in the representation $\rho_{p,\xi,\eta}$ we obtain the property $\phi(\lambda)=\phi(\lambda+1)$. The only polynomial in $\lambda$ with these properties is $\phi(\lambda)=\lambda^p-\lambda-\xi^p\eta^p$, which proves (\ref{detf}).

  Let $\mu_1,...,\mu_p\in\overline{\Z/(p)(\xi,\eta)}$ be the eigenvalues of $\rho_{p,\xi,\eta}(x\partial_x)$. It follows from (\ref{detf}) that
  $$(\lambda-\mu_1)...(\lambda-\mu_p)=\lambda^p-\lambda-\xi^p\eta^p.$$
  Taking $j$-th derivative with respect to $\lambda$ of logarithm of the l.h.s. and r.h.s. of this identity we obtain 
  $$\frac{1}{(\lambda-\mu_1)^j}+...+\frac{1}{(\lambda-\mu_p)^j}=\frac{1}{(\xi^p\eta^p+\lambda-\lambda^p)^j},~~~j=1,...,p$$
  which gives (\ref{trf2}) for $i=0$.

  Equation (\ref{trf1}) in the case $i=0$ can be proved similarly: we need to express the sums of powers $\mu_1^j+...+\mu_p^j$ in terms of elementary symmetric functions and use formula (\ref{detf}) for the characteristic polynomial of $\rho_{p,\xi,\eta}(x\partial_x)$.

  It follows from (\ref{comrel}) that 
  $$(x\partial_x-\lambda)\cdot x^i(x\partial_x-\lambda)^{-j+1}\cdot (x\partial_x-\lambda)^{-1}=x^i(x\partial_x-\lambda)^{-j+1}\cdot (x\partial_x-\lambda+i)(x\partial_x-\lambda)^{-1}=$$
  $$=x^i(x\partial_x-\lambda)^{-j+1}+ix^i(x\partial_x-\lambda)^{-j}.$$
  Taking trace in the representation $\rho_{p,\xi,\eta}$ gives
  $$i~\text{tr}~(\rho_{p,\xi,\eta}(x^i(x\partial_x-\lambda)^{-j}))=0$$
  which gives the equations (\ref{trf1}) and (\ref{trf2}) in the case $i\ne 0$ mod $p$.

  $\square$

  Using this Lemma and the expression (\ref{A}) we obtain
  $$\text{tr}~(\rho_{p,\xi,\eta}(A))=\sum_{j=-d_1-...-d_N}^{d_1+...+d_N}\text{tr}~(\rho_{p,\xi,\eta}(x^jg_j(x\partial_x)))=\text{tr}~(\rho_{p,\xi,\eta}(g_0(x\partial_x)))=T_{\varepsilon}g_0(\varepsilon) |_{\varepsilon=-\xi^p\eta^p}$$
  where we assume $p>d_1+...+d_N$ and use (\ref{tr}) and the definition of $T_{\varepsilon}$.

  On the other hand, from (\ref{A}) it follows that $\rho_{\varepsilon}(A)x^{\varepsilon}=\sum_{j=-d_1-...-d_N}^{d_1+...d_N}g_j(\varepsilon)x^{\varepsilon+j}$ and therefore $\text{Coef}_{x^{\varepsilon}}(\rho_{\varepsilon}(A)x^{\varepsilon})=g_0(\varepsilon)$ which gives the statement of the theorem.
  
  $\square$

  {\bf Theorem 2.2.} Let $D=P+Q$ where $P$ is given by (\ref{P}) and $Q\in t\cdot \text{Diff}[[t]]$, where $t$ is a formal parameter. In other words, 
  $$Q=\sum_{j=1}^{\infty}Q_jt^j,~~~Q_j\in\text{Diff}.$$
  Then for any $N\geqslant 1$ the exists $C_N>0$ such that for all prime numbers $p>C_N$ we have 
  \begin{equation}\label{detD}\det(\rho_{p,\xi,\eta}(D))=(-1)^k\varepsilon^k e^{T_{\varepsilon}\text{Coef}_{x^{\varepsilon}}(\rho_{\varepsilon}(\log(1+P^{-1}Q))x^{\varepsilon})}|_{\varepsilon=-\xi^p\eta^p}~~~\text{mod}~p,~t^N.
  \end{equation}
   We understand $\log$ as a power series and $T_{\varepsilon}$ is applied to coefficients of all powers of $t$.

  {\bf Proof.}
  Notice that $\det(\rho_{p,\xi,\eta}(P))=\xi^{pk}\eta^{pk},$ and we can write
  $$\det(\rho_{p,\xi,\eta}(D))=\det(\rho_{p,\xi,\eta}(P(1+P^{-1}Q)))=\xi^{kp}\eta^{kp}\det(\rho_{p,\xi,\eta}(1+P^{-1}Q))=$$
  $$=\xi^{kp}\eta^{kp}e^{\text{tr}~(\log_{(p)}(\rho_{p,\xi,\eta}(1+P^{-1}Q)))}~~~\text{mod}~t^{p-1}$$
  where $\log_{(p)}$ here is a truncated $\log$ and understood as a power series $\log_{(p)}(1-x)=-x-\frac{x^2}{2}-...-\frac{x^{p-1}}{p-1}$. Using Theorem 2.1 to calculate the trace, we get the formula (\ref{detD}).

  $\square$

  {\bf Remark 2.1.} The above formulas can be generalized as follows. Let $A,~P$ be defined as in Theorem 2.1, but assume that $r_1,...,r_k$ are formal parameters.
  Define a linear operator $T_{p,\varepsilon}:~{\bf k}(r_1,...,r_k)[\varepsilon,(\varepsilon-r_j+m)^{-1},m\in\Z,1\leqslant j\leqslant k]\to {\bf k}(r_1,...,r_k)[\varepsilon,(\varepsilon-r_j+r_j^p)^{-1},1\leqslant j\leqslant k]$ by
  $$T_{p,\varepsilon}(\varepsilon^l)=0,~l\geqslant 0,~~~~~~~T_{p,\varepsilon}\Bigg(\frac{1}{(\varepsilon-r_i+m)^j}\Bigg)=\frac{1}{(\varepsilon-r_i+r_i^p)^j},~j>0,1\leqslant i\leqslant k,m\in\Z.$$
  Then the following formulas hold for all large enough primes $p$
  $$\text{tr}~ (\rho_{p,\xi,\eta}(A))=T_{p,\varepsilon} \text{Coef}_{x^{\varepsilon}}(\rho_{\varepsilon}(A)x^{\varepsilon}) |_{\varepsilon=-\xi^p\eta^p}~~~\text{mod}~p$$
  $$\det(\rho_{p,\xi,\eta}(D))=(-1)^k\prod_{i=1}^k(\varepsilon-r_i+r_i^p) \cdot e^{T_{p,\varepsilon}\text{Coef}_{x^{\varepsilon}}(\rho_{\varepsilon}(\log(1+P^{-1}Q))x^{\varepsilon})}|_{\varepsilon=-\xi^p\eta^p}~~~\text{mod}~p,~t^N~~~\text{for}~p>C_N.$$
  Proofs are similar to ones of Theorems 2.1 and 2.2, using formulas (\ref{detf}), (\ref{trf1}), (\ref{trf2}).

  \section{Formal solutions in characteristic zero and another explicit formula for determinants}

  Let ${\bf k}$ be a field of characteristic zero, ${\bf \bar{k}}$ its algebraic closure, 
  $$P=(x\partial_x-r_1)...(x\partial_x-r_k),~~~~~~~r_1,...,r_k\in{\bf k}$$
  be a monic polynomial in $x\partial_x$ on ${\bf k}$, and \begin{equation}\label{Q}
  Q=\sum_{i=1}^{\infty}Q_it^i,~~~~Q_i\in{\bf k}(\!(x)\!)[\partial_x]
  \end{equation}
  where 
  \begin{equation}\label{Qi}
  Q_i=\sum_{j=-N_i}^{\infty}\sum_{l=0}^{M_i} q_{i,j,l}x^j\partial_x^l,~~~~q_{i,j,l}\in {\bf k}.
  \end{equation}
  Notice that $N_i,M_i$ are not bounded and can increase with $i$.

  Define an equivalence relation $r_i \sim r_j$ iff $r_i-r_j\in\Z$. Let $s_1,...,s_l$ be the smallest representatives in each class of equivalence\footnote{In other words, $s_i-s_j\notin\Z$ for $i\ne j$ and each $r_i$ has a form $r_i=s_j+m$ for some $j=1,...,l$ and $m\geqslant 0$.}.  We can write
  \begin{equation}\label{pol1}
  P=\prod_{i=1}^l\prod_{j\geqslant 0} (x\partial_x-s_i-j)^{a_{i,j}}
  \end{equation}
  where $a_{i,j}\geqslant 0$.  Let $k_i=\sum_{j\geqslant 0} a_{i,j}$. Then $k_1+...+k_l=\sum_{i,j}a_{i,j}=k$.

  {\bf Theorem 3.1.} Let 
  $$V_i=x^{s_i}{\bf k}(\!(x)\!)[\log x][[t]],~~~i=1,...,l.$$
  Then the space of solutions of the equation
  \begin{equation}\label{eq}
  (P+Q)f=0,~~~f\in x^{r}{\bf k}(\!(x)\!)[\log x][[t]],~~~ r\in {\bf k}
  \end{equation}
  is nonzero iff $r-s_i\in\Z$ for some $i=1,...,l$, so in this case, without loss of generality, we can set $r=s_i$ and therefore $f\in V_i$. The space of solutions of equation (\ref{eq}) with $f\in V_i$ is a free $k_i$-dimensional module over ${\bf k}[[t]]$. Moreover, this module after an extension of scalars has a basis $f_{i,1},...,f_{i,k_i}$ such that 
  \begin{equation}\label{fij}f_{i,j}=x^{s_i+\lambda_{i,j}}g_{i,j}\in x^{s_i+\lambda_{i,j}}{\bf \bar{k}}(\!(x)\!)[[t^{1/L}]][\log x],~~~\lambda_{i,j}\in t^{1/L}\cdot{\bf \bar{k}}[[t^{1/L}]]
  \end{equation}
  for some $L\geqslant 1$, where degrees of elements $g_{i,j}\in {\bf \bar{k}}(\!(x)\!)[[t^{1/L}]][\log x]$ as polynomials in $\log x$ are smaller than $k_i$, and $x^{s_i+\lambda_{i,j}}$ is understood as
  $$x^{s_i+\lambda_{i,j}}=x^{s_i}\cdot e^{\lambda_{i,j}~\log x}=x^{s_i}\sum_{m=0}^{\infty}\frac{\lambda_{i,j}^m(\log x)^m}{m!}\in x^{s_i}{\bf \bar{k}}(\!(x)\!)[\log x][[t^{1/L}]].$$

  {\bf Proof.} Let $f=f_0+f_1t+f_2t^2+...$ be a nonzero solution of (\ref{eq}) where $f_0,f_1,f_2,...\in x^{r}{\bf k}(\!(x)\!)[\log x]$. Without loss of generality we can assume\footnote{If $f=f_mt^m+f_{m+1}t^{m+1}+...$ is a solution with $f_m\ne 0$,  then $f=f_m+f_{m+1}t+...$ is also a solution.} $f_0\ne 0$. Equation (\ref{eq}) is equivalent to the system:
  \begin{equation}\label{eqf}Pf_m=-Q_1f_{m-1}-Q_2f_{m-2}-...-Q_mf_0,~~~m=0,1,2,...,
  \end{equation}
  see (\ref{Q}). In particular, for $m=0$ we have the equation $Pf_0=0$ which can be solved explicitly. This equation has nonzero solutions iff $r-s_i\in\Z$ for some $i=1,...,l$. So we can set $r=s_i$ without loss of generality. Moreover, one can check that the linear operator 
  $$P:~~~x^{s_i}{\bf k}(\!(x)\!)[\log x]~\to~x^{s_i}{\bf k}(\!(x)\!)[\log x]$$
  has $k_i$-dimensional kernel with a basis $\{x^{s_i+j}(\log x)^s,~j\geqslant 0,~0\leqslant s\leqslant a_{i,j}-1\}$ and acts epimorphically. In other words, we have
  \begin{equation}\label{opP}
  \text{Ker} P=\text{span}(x^{s_i+j}(\log x)^s,~j\geqslant 0,~0\leqslant s\leqslant a_{i,j}-1),~~~~~~~\text{Im} P=x^{s_i}{\bf k}(\!(x)\!)[\log x].
  \end{equation}
  So we can choose  $f_0\in\text{Ker} P$
  arbitrarily and for each $m=1,2,...$ we can solve (\ref{eqf}) for $f_m$ assuming that $f_0,...,f_{m-1}$ are already found. Moreover this solution $f_m$ is unique up to an arbitrary element in $\text{Ker} P$. Fix a direct sum decomposition\footnote{For example, one can set $W=\{x^{s_i}\sum_{j,s}q_{j,s}x^j(\log x)^s;~~~q_{j,s}=0 \text{~if~} j\geqslant 0,~0\leqslant s\leqslant a_{i,j}-1\}$.} $x^{s_i}{\bf k}(\!(x)\!)[\log x]=\text{ker} P\oplus W$. It follows from (\ref{opP}) that for each $j\geqslant 0,~s=0,...,a_{i,j}-1$ such that $a_{i,j}>0$ there exist a unique solution $f=h_{i,j,s}$ such that
  $f_0=x^{s_i+j}(\log x)^s$ and $f_m\in W$ for $m\geqslant 1$. Moreover, any solution $f\in V_i$ can be written uniquely in the form $f=\sum_{j,s}b_{j,s}g_{j,s}$ for some $b_{j,s}\in {\bf k}[[t]]$.

  Let $D_{\log}$ be a derivation of the algebra 
  $$\bigoplus_{r\in{\bf k}}x^{r}{\bf k}(\!(x)\!)[\log x][[t]] $$
  defined by $D_{\log}(x^r)=D_{\log}(x)=D_{\log}(t)=0$, $D_{\log}(\log x)=1$. 
  Notice that operators $P+Q,~D_{\log}$ both act on $V_i$ and $[P+Q,D_{\log}]=0$ because the coefficients of $P+Q$ do not contain $\log x$. Therefore, the operator $D_{\log}$ acts on the space of solutions of (\ref{eq}) in $V_i$ for all $i=1,...,l$. Let $g_{i,1},...,g_{i,k_i}$ be a ${\bf k}[[t]]$-basis of this space of solutions. We have
  $$D_{\log}g_{i,j}=\sum_{m=1}^{k_i}c_{i,j,m}g_{i,m},~~~j=1,...,k_i$$
  for some $c_{i,j,m}\in{\bf k}[[t]]$. The characteristic polynomial of the $k_i\times k_i$ matrix $A_i=(c_{i,j,m})_{j,m}$ can be factorized into linear factors after a certain extension of scalars. In fact, $A_i|_{t=0}$ is nilpotent because it is a matrix of the operator $D_{\log}$ acting on $\text{Ker}P$. So its characteristic polynomial $\det(A_i-\lambda)|_{t=0}=(-1)^{k_i}\lambda^{k_i}$. Therefore, for some $L\in \N$ such that $1\leqslant L\leqslant k_i$ we have
  $$\det(A_i-\lambda)=(-1)^{k_i}(\lambda-\lambda_{i,1})...(\lambda-\lambda_{i,k_i})$$
  where $\lambda_{i,j}\in t^{1/L}\cdot{\bf \bar{k}}[[t^{1/L}]]$. After such extension of scalars we can choose a basis in the space of solutions of (\ref{eq}) in the vector space $V_i\otimes t^{1/L}\cdot{\bf \bar{k}}(\!(t^{1/L})\!)$ consisting of generalized eigenvectors with eigenvalues $\lambda_{i,j}$. Notice that if $f\in x^{s_i}{\bf \bar{k}}(\!(x)\!)[\log x][[t^{1/L}]]$ is a generalized eigenvector, $(D_{\log}-\lambda_{i,j})^Mf=0$, then $f\in x^{s_i+\lambda_{i,j}}{\bf \bar{k}}(\!(x)\!)[[t^{1/L}]][\log x]$ and its degree as a polynomial in $\log x$ is less than $M$.

  $\square$

  We define a linear operator $T_{\varepsilon}:~{\bf k}(\varepsilon)\to \varepsilon^{-1}{\bf k}[\varepsilon^{-1}]$ by formulas (\ref{T}). Notice that in order to apply $T_{\varepsilon}$ to a rational function $f(\varepsilon)\in {\bf k}(\varepsilon)$ we need to embed ${\bf k}\subset \overline{{\bf k}}$ in its algebraic closure, factorize the denominator of $f(\varepsilon)$, use the partial fraction decomposition of $f(\varepsilon)$ and apply the formulas (\ref{T}). However, the operator $T_{\varepsilon}$ is invariant with respect to the Galois group $Gal(\overline{{\bf k}}/{\bf k})$, so we have $T_{\varepsilon}(f(\varepsilon))\in {\bf k}(\varepsilon)$. 

  {\bf Theorem 3.2.} We have 
  \begin{equation}\label{pol}\varepsilon^ke^{T_{\varepsilon}(\text{Coef}_{x^{\varepsilon}}(\log(1+P^{-1}Q)(x^{\varepsilon})))}=\prod_{i=1}^l\prod_{j=1}^{k_i} (\varepsilon-\lambda_{i,j})
  \end{equation}
  where $\lambda_{i,j}\in t^{1/L}\cdot{\bf \bar{k}}[[t^{1/L]}]$ are defined in Theorem 3.1, and in the l.h.s. we understand $\log$ as power series which can be written as a power series in $t$, and $T_{\varepsilon}$ is applied to each coefficient of this power series. In particular, the expression in the l.h.s. of (\ref{pol}) is a monic polynomial in $\varepsilon$ of degree $k$.

  {\bf Corollary 3.1.} Let $s_1,...,s_l\in\Q$ and $s_i-s_j\notin\Z$ for $i\ne j$. Then combining the Theorems 2.2, 3.1, 3.2 we get
$$\det(\rho_{p,\xi,\eta}(D))=(-1)^kw_1(-\xi^p\eta^p)...w_l(-\xi^p\eta^p)~~~\text{mod}~p,~t^N$$
for all primes $p>C_N$.

More generally, if $s_1,...,s_l\in{\bf k}$ and $s_i-s_j\notin\Z$ for $i\ne j$, then
$$\det(\rho_{p,\xi,\eta}(D))=(-1)^kw_1(-\xi^p\eta^p-s_1+s_1^p)...w_l(-\xi^p\eta^p-s_l+s_l^p)~~~\text{mod}~p,~t^N$$
for all primes $p>C_N$.

  {\bf Proof.} It is sufficient to prove the theorem under the following assumptions:
  
  {\bf 1.} The field ${\bf k}$ is algebraically closed, and $L=1$. 

  {\bf 2.} For any given $i$, all $\lambda_{i,j}$ are distinct.

  {\bf 3.} $Q_j\in {\bf k}[x^{\pm 1},\partial_x],~j=1,2,3,...$

  Indeed, for fixed $P$ the statement of the theorem is a collection of certain polynomial identities that involve coefficients of operators $Q_1,Q_2,...$. In each identity only finitely many coefficients appear. Therefore, it is sufficient to prove the theorem under the assumption {\bf 3}. All $\lambda_{i,j}$ are distinct if the coefficients of $Q_1,Q_2,...$ are in a generic position, which gives the assumption {\bf 2.} Finally, the first assumption can be achieved by extension of the field ${\bf k}$ and change of variable $t\mapsto t^L$.

  It follows from these assumptions that 
  \begin{equation}\label{ass}
  \lambda_{i,j}\in t{\bf k}[[t]],~~~f_{i,j}\in x^{s_i+\lambda_{i,j}}{\bf k}[x,x^{-1}][[t]],~~~(P+Q)f_{i,j}=0
  \end{equation}
  for all $\lambda_{i,j},~f_{i,j}$ because the matrices $A_i$ of the operator $P+Q$ acting in the space $V_i$ are diagonal in the basis $f_{i,1},...,f_{i,k_i}$ for all $i=1,...,l$.

  Fix an arbitrary, large enough integer $M$ such that all $\lambda_{i,j}$ are distinct modulo $t^M$ and work over ${\bf k}[[t]]/(t^M)={\bf k}[t]/(t^M)$. After that we have
  \begin{equation}\label{fij1}\lambda_{i,j}\in t{\bf k}[t]/(t^M),~~~~~~~f_{i,j}=x^{s_i+\lambda_{i,j}}g_{i,j},~~~g_{i,j}\in{\bf k}[x,x^{-1}][t]/(t^M).
  \end{equation}
  where $g_{i,j}$ are Laurent polynomials in $x$ because $Q_j$ are polynomials in $x$. 

  The operator $P+Q$ acts linearly on $W_{\varepsilon}$, which is a free module over ${\bf k}[t]/(t^M)$ with the basis $\{v_n=x^{\varepsilon+n};~n\in\Z\}$ as a linear operator $\rho_{\varepsilon}(P+Q)$, where $\varepsilon$ is a formal parameter, see (\ref{repep}). Let $B(\varepsilon)$ be the corresponding matrix of infinite size.  The matrix elements of $B(\varepsilon)$ are polynomials in $\varepsilon$. Let us identify solutions $f_{i,j}$ written in the form (\ref{fij1}) with elements $e_{i,j}\in W_{\varepsilon}$ as follows: if $g_{i,j}=\sum_nc_{n,i,j}x^n$, then the corresponding element $e_{i,j}\in W_{\varepsilon}$ reads $e_{i,j}=\sum_nc_{n,i,j}v_n$. It follows from (\ref{ass}) that $B(s_i+\lambda_{i,j})e_{i,j}=0$ for all $e_{i,j}$.

  Notice that the matrix elements of $B(\varepsilon)$ vanish outside of finitely many shifts of the main diagonal, i.e. we have $B(\varepsilon)_{m,n}=0$ if $|n-m|>C$ where $C$ depends on $M$. Denote by $B_N(\varepsilon)$ for $N\geqslant 1$ the $(2N+1)\times (2N+1)$ submatrix of $B(\varepsilon)$ which corresponds to the subset $\{v_n;~-N\leqslant n\leqslant N\}\subset \{v_n;~n\in\Z\}$. We have $B_N(s_i+\lambda_{i,j})e_{i,j}=0$ for large enough $N$. Therefore, we have \begin{equation}\label{rootsB}\det B_N(s_i+\lambda_{i,j})=0,~~~N\gg 1.
  \end{equation}

  Let us investigate the structure of $\det B_N(\varepsilon)$ for $N\gg 1$.

  {\bf Proposition 3.1.} There exist $a\geqslant 0$ independent of $N$, such that if $N\geqslant a$, then we have 
 \begin{equation}\label{detb}\det B_{N}(\varepsilon)=(1+\phi_1(\varepsilon-N))(1+\phi_2(\varepsilon+N))f_1(\varepsilon-N)\cdot \prod_{n=a-N}^{N-a}\prod_{i=1}^lw_i(\varepsilon-s_i+n)\cdot f_2(\varepsilon+N)
 \end{equation}
 where $f_1(\varepsilon),f_2(\varepsilon),w_i(\varepsilon),\phi_1(\varepsilon),\phi_2(\varepsilon)\in {\bf k}[t]/(t^M)[\varepsilon]$ are such that

 \begin{align}1)&~~~ f_1(\varepsilon)=f_1^{(0)}(\varepsilon)\text{~~~mod~}t\\
2)&~~~ f_2(\varepsilon)=f_2^{(0)}(\varepsilon)\text{~~~mod~}t\\
 3)&~~~w_i(\varepsilon)=\varepsilon^{k_i}\text{~~~mod~}t,~~~\text{and are monic polynomials}\\
 4)&~~~\phi_1(\varepsilon),~\phi_2(\varepsilon)=0\text{~~~mod~}t\end{align}
and  $f_1(\varepsilon),f_2(\varepsilon)\in{\bf k}[\varepsilon]$ are monic polynomials  which are uniquely defined by the equality
\begin{equation}\label{f12}
\prod_{n=-N}^N\prod_{i,j}(\varepsilon-s_i-j+n)^{a_{i,j}}=f_1^0(\varepsilon-N)\cdot\prod_{n=a-N}^{N-a}\prod_{i=1}^l(\varepsilon-s_i+n)^{k_i}\cdot f_2^0(\varepsilon+N).
\end{equation}
Notice that the roots of $f_1^0(\varepsilon),f_2^0(\varepsilon)$ belong to $\bigcup_i(s_i+\Z)\subset{\bf k}$.

 {\bf Proof.} Recall that the operator $Q$ vanishes at $t=0$, so the matrix $B(\varepsilon)$ at $t=0$ is diagonal because it corresponds to the operator $P$. We have 
 \begin{equation}\label{detb0}
 B(\varepsilon)_{m,n}=\delta_{m,n}~\prod_{i,j}(\varepsilon-s_i-j+n)^{a_{i,j}}~~~\text{mod}~t.
 \end{equation}
 Therefore, for determinant of its submatrix we have
 $$\det B_N(\varepsilon)=\prod_{n=-N}^N\prod_{i,j}(\varepsilon-s_i-j+n)^{a_{i,j}}~~~\text{mod}~t.$$
 Using the equation (\ref{f12}) we obtain
 \begin{equation}\label{b0}
 \det B_N(\varepsilon)=f_1^0(\varepsilon-N)\cdot\prod_{n=a-N}^{N-a}\prod_{i=1}^l(\varepsilon-s_i+n)^{k_i}\cdot f_2^0(\varepsilon+N)~~~\text{mod}~t.
 \end{equation}
 This means that equation (\ref{detb}) holds modulo $t$. 

 Notice that equation (\ref{detb0}) gives a decomposition of $\det B_N(\varepsilon)$ at $t=0$ into a product of $2(N-a+1)+2$ polynomials without common roots. We will apply to this decomposition the following lemma.

{\bf Lemma 3.1.} Let $\bar{g}_1,...\bar{g}_r\in {\bf k}[\varepsilon]$ be monic polynomials without common roots. Let $g\in {\bf k}[t]/(t^M)[\varepsilon]$ be a polynomial such that $g=\prod_{i=1}^r\bar{g}_i\text{~mod~}t$. Then there exist unique monic polynomials $g_i\in {\bf k}[t]/(t^M)[\varepsilon]$, $g_i=\bar{g}_i\text{~mod~}t$ and $h\in {\bf k}[t]/(t^M)[\varepsilon]$, $h=1\text{~mod~}t$ such that $g=h\prod_{i=1}^rg_i$.

{\bf Proof.} Assume that we have found such a decomposition modulo $t^l$ for some $l\geqslant 1$ and let us find it modulo $t^{l+1}$. Choose some lifts $\tilde{g}_1,...\tilde{g}_r,\tilde{c}\in {\bf k}[t]/(t^M)[\varepsilon]/t^{l+1}$ of $g_1,...,g_r,c\in {\bf k}[t]/(t^M)[\varepsilon]/t^l$, such that $\tilde{g}_1,...\tilde{g}_r$ are monic polynomials in $\varepsilon$. We have 
\begin{equation} \label{exp}
\tilde{c}\prod_{i=1}^r\tilde{g}_i-g\in t^{l+1}{\bf k}[t]/(t^M)[\varepsilon].
\end{equation}
We want to correct our lift so that the expression (\ref{exp}) becomes zero modulo $t^{l+2}$. After replacing $\tilde{c}\mapsto \tilde{c}+\delta\tilde{c},~\tilde{q}_i\mapsto \tilde{g}_i+\delta\tilde{g}_i$ the expression (\ref{exp}) modulo $t^{l+2}$ reads
\begin{equation}\label{expn}
\tilde{c}\prod_{i=1}^r\tilde{g}_i-g+\delta\tilde{c}~\prod_{i=1}^r\tilde{g}_i+\sum_{i=1}^r\tilde{c}~\delta\tilde{g}_i\prod_{j\ne i}\tilde{g}_j\in I^{l+1}[\varepsilon]/(t^{l+2}).
\end{equation}
Here $\delta\tilde{c},\delta\tilde{g}_i\in t^{l+1}{\bf k}[t]/(t^M)[\varepsilon]$ and $\deg\delta\tilde{g}_i<\deg\bar{g}_i$ because $\tilde{g}_i+\delta\tilde{g}_i$ should be monic polynomials as well as $\tilde{g}_i$.
The assumption that $\bar{g}_i,\bar{g}_j$ are coprime polynomials implies (by chinese remainder theorem) that
$${\bf k}[\varepsilon]/\prod_{i=1}^r\bar{q}_i \cong \prod_{i=1}^r{\bf k}[\varepsilon]/(\bar{g}_i). $$
It follows that $\{\varepsilon^{a_i}\prod_{j\ne i}\bar{g}_i;~1\leqslant i\leqslant r,0\leqslant a_i\leqslant\deg\bar{g}_i-1\}$ form a ${\bf k}$-basis of ${\bf k}[\varepsilon]/\prod_{i=1}^r\bar{q}_i$. Hence any element $p\in {\bf k}[\varepsilon]$ can be written uniquely as 
\begin{equation}\label{comb}
p=\sum_{i=1}^r\sum_{a_i=0}^{\deg \bar{g}_i-1}\phi_{i,a_i}\varepsilon^{a_i}\prod_{j\ne i}\bar{g}_i +\prod_{i=1}^r\bar{g}_i\sum_{a=0}^{\infty}\psi_a\varepsilon^a
\end{equation} 
where $\phi_{i,a_i},\psi_a\in {\bf k}$ and $\psi_a=0$ for $a\gg 0$.
We can set $p=\tilde{c}\prod_{i=1}^r\tilde{g}_i-g\in t^{l+1}{\bf k}[t]/(t^M)[\varepsilon]/(t^{l+2})$ in (\ref{comb}) and get the unique modulo $t^{l+2}$ correction $\delta\tilde{c},\delta\tilde{g}_i$ as coefficients of expansion (\ref{comb}). 

$\square$

{\bf Remark 3.1.} In fact, a more precise statement is true. Let $R\subset {\bf k}$ be a subring containing coefficients of polynomials $g,\bar{g}_1,...,\bar{g}_r$. Assume that the resultants $\text{Res}(\bar{g}_i,\bar{g}_j),~i<j$ are invertible in $R$. Then all the coefficients of the polynomials $g_1,...,g_r$ belong to $R$. 
 
Applying Lemma 3.1 to (\ref{b0}) we get the determinant $\det B_{N}(\varepsilon)$ in the form
\begin{equation}\label{b1}
\det B_{N}(\varepsilon)=h_{N}(\varepsilon)f_{1,N}(\varepsilon-N)\cdot \prod_{n=a-N}^{N-a}\prod_{i=1}^lw_{i,N,n}(\varepsilon-s_i+n)\cdot f_{2,N}(\varepsilon+N)
\end{equation}
where $f_{1,N}(\varepsilon),w_{i,N,n}(\varepsilon),f_{2,N}(\varepsilon)$ are monic polynomials in ${\bf k}[t]/(t^M)[\varepsilon]$ such that
$$f_{1,N}(\varepsilon)=f_1^0(\varepsilon)~~~\text{mod}~t,~~~~~~~f_{2,N}(\varepsilon)=f_2^0(\varepsilon)~~~\text{mod}~t,~~~~~~~w_{i,N,n}(\varepsilon)=\varepsilon^{k_i}~~\text{mod}~t,$$
and $h_{N}(\varepsilon)\in {\bf k}[t]/(t^M)[\varepsilon]$ is equal to 1 modulo $t$.

Introduce the notation $B_N^{(0)}(\varepsilon)=B_N(\varepsilon)|_{t=0}$. Recall that $B_N^{(0)}(\varepsilon)$ is a diagonal matrix, it corresponds to the operator $P$ acting in the free module over ${\bf k}[t]/(t^M)$ with the basis $\{x^{\varepsilon+n};~n=-N,...,N\}$. It follows from the identity $P+Q=P(1+P^{-1}Q)$ that 
$$\det B_N(\varepsilon)=\det B_N^{(0)}(\varepsilon)\det (1+A_N(\varepsilon)$$
where $A_N(\varepsilon)=B_N^{(0)}(\varepsilon)^{-1}B_N(\varepsilon)-1$. Therefore, using (\ref{b0}) we obtain
\begin{equation}\label{dbf} \det B_N(\varepsilon)=f_1^0(\varepsilon-N)\cdot\prod_{n=a-N}^{N-a}\prod_{i=1}^l(\varepsilon-s_i+n)^{k_i}\cdot f_2^0(\varepsilon+N)\det(1+A_N(\varepsilon)).
\end{equation}
Matrix elements of $A_N(\varepsilon)$ belong to ${\bf k}[\varepsilon,~\frac{1}{\varepsilon-s_i+m},i=1,...,l,~m\in\Z]~t{\bf k}[t]/(t^M)$ and has the form
$$A_N(\varepsilon)_{m,n}=a_{n-m}(\varepsilon+n)$$
where $a_n(\varepsilon)=0$ for $n>C$ and $C$ is independent of $N$.

{\bf Lemma 3.2.} Let $A(\varepsilon)$ be an arbitrary infinite matrix such that its matrix elements $A(\varepsilon)_{m,n},~m,n\in\Z$ have the form $A(\varepsilon)_{m,n}=a_{n-m}(\varepsilon+n)$ where $a_n(\varepsilon)\in {\bf k}[\varepsilon,~\frac{1}{\varepsilon-s_i+m},i=1,...,l,~m\in\Z]~t{\bf k}[t]/(t^M)$ and $a_n(\varepsilon)=0$ for $|n|\gg 0$. Then for finite submatrices $A_N(\varepsilon)$ corresponding to indices $n\in\{-N,...,N\}$ we have 
$\log\det(1+A_N(\varepsilon))\in {\bf k}[\varepsilon,~\frac{1}{\varepsilon-s_i+m},i=1,...,l,~m\in\Z]~t{\bf k}[t]/(t^M)$ and for $N\gg 0$
\begin{equation}\label{det}
\log\det(1+A_N(\varepsilon))=h_1(\varepsilon-N)+\sum_{n=a-N}^{N-a}\Big(\sum_{r=1}^{r_0}\sum_{i=1}^l\frac{c_{i,r}}{(\varepsilon-s_i+n)^r}+\sum_{\tilde{r}=0}^{\tilde{r}_0}\tilde{c}_{\tilde{r}}(\varepsilon+n)^{\tilde{r}}\Big)+h_2(\varepsilon+N)
\end{equation}
for some $a,b,r_0,\tilde{r}_0\geqslant 0$, where $c_{i,r},\tilde{c}_r\in t{\bf k}[t]/(t^M)$ and $h_1(\varepsilon),h_2(\varepsilon)\in {\bf k}[\varepsilon,~\frac{1}{\varepsilon-s_i+m},i=1,...,l,~m\in\Z]~t{\bf k}[t]/(t^M)$.

{\bf Proof.} We have $\log\det(1+A_N(\varepsilon))=\text{tr}~\log(1+A_N(\varepsilon))=$
$$\sum_{l=1}^{\infty}\frac{(-1)^{l+1}}{l}\sum_{-N\leqslant n_1,...,n_l\leqslant N}A_{n_1,n_2}(\varepsilon)A_{n_2,n_3}(\varepsilon)...A_{n_{l-1},n_l}(\varepsilon)A_{n_l,n_1}(\varepsilon)=$$
\begin{equation}\label{log}
\sum_{l=1}^{\infty}\frac{(-1)^{l+1}}{l}\sum_{-N\leqslant n_1,...,n_l\leqslant N}a_{n_2-n_1}(\varepsilon+n_1)...a_{n_l-n_{l-1}}(\varepsilon+n_{l-1})a_{n_1-n_l}(\varepsilon+n_l).
\end{equation}
Notice that all terms in the r.h.s. of (\ref{log}) with $l>M$ vanish because $a_n(\varepsilon)=0$ mod $t$ and we work modulo $(t^M)$, and  $a_n(\varepsilon)=0$ for $|n|\gg 0$. Therefore, the summation in $l$ is finite and only a finite number of sequences of indexes $(n_2-n_1,n_3-n_2,...,n_l-n_{l-1},n_1-n_l)$ appears in (\ref{log}) and this set of sequences of indexes does not depend on $N$ for $N\gg 0$. Using partial fraction decomposition we can rewrite (\ref{log}) as a finite $t{\bf k}[t]/(t^M)$-linear combination of sums of the form
$$\sum_{n=\alpha-N}^{N-\beta}\frac{1}{(\varepsilon-s_i+n)^r},~i=1,...,l,~~~1\leqslant r \text{~~~and~~~}\sum_{n=\alpha-N}^{N-\beta}(\varepsilon+n)^{\tilde{r}},~\tilde{r}\geqslant 0$$
where parameters $\alpha,\beta\geqslant 0$ vary.

Let $a$ be the maximum of values of $\alpha$ and $\beta$. Then the contribution of terms $\frac{1}{(\varepsilon-s_i+n)^r},~(\varepsilon+n)^r$ for fixed non-zero $r$ and $a-N\leqslant n\leqslant N-a$ gives the middle term in (\ref{det}), while the combination of these terms with $n<a-N$ (resp. $n>N-a$) gives $h_1(\varepsilon-N)$ (resp. $h_2(\varepsilon+N)$). 

$\square$

Let us apply Lemma 3.2 to the infinite matrix $A(\varepsilon)=B^{(0)}(\varepsilon)^{-1}B(\varepsilon)-1$ where $B^{(0)}(\varepsilon)$ is the diagonal matrix of the operator $P$. Using (\ref{dbf}), we get
\begin{equation}\label{lar}
\det B_{N}(\varepsilon)=f_1^0(\varepsilon-N)\cdot\prod_{n=a-N}^{N-a}\prod_{i=1}^l(\varepsilon-s_i+n)^{k_i}\cdot f_2^0(\varepsilon+N)\times
\end{equation}
$$\exp(h_1(\varepsilon-N))\exp(h_2(\varepsilon+N))\prod_{n=N-a}^{N+b}\exp\Big(\sum_{r=1}^{r_0}\sum_{i=1}^l\frac{c_{i,r}}{(\varepsilon-s_i+n)^r}+\sum_{\tilde{r}=0}^{\tilde{r}_0}\tilde{c}_{\tilde{r}}(\varepsilon+n)^{\tilde{r}}\Big).$$
Let us consider Laurent expansion of $\det B_{N}(\varepsilon)$ at $\varepsilon=s_i-n$ where $N\pm n\gg 0$. The formula (\ref{lar}) gives $(\varepsilon-s_i+n)^{k_i}\exp\Big(\sum_{r=1}^{r_0}\frac{c_{i,r}}{(\varepsilon-s_i+n)^r}\Big)$ times an invertible element of ${\bf k}[t]/(t^M)[[\varepsilon-s_i+n]]$, where $c_{i,r}\in t{\bf k}[t]/(t^M)$. Using the fact that $\det B_{N}(\varepsilon)$ is a polynomial in $\varepsilon$ and hence has no poles at $\varepsilon=s_i-n$ we conclude that there exists unique monic polynomial $w_i(\varepsilon)\in {\bf k}[t]/(t^M)[\varepsilon]$ of degree $k_i$ such that 
$$(\varepsilon-s_i+n)^{k_i}\exp\Big(\sum_{r=1}^{r_0}\frac{c_{i,r}}{(\varepsilon-s_i+n)^r}\Big)=w_i(\varepsilon-s_i+n).$$
On the other hand, recall another factorization of $\det B_N(\varepsilon)$ given by (\ref{b1}). Therefore, we get an equality
$$f_1^0(\varepsilon-N)\exp(h_1(\varepsilon-N))\cdot f_2^0(\varepsilon+N)\exp(h_2(\varepsilon+N))\cdot \prod_{n=a-N}^{N-a}\exp(\sum_{\tilde{r}_0=0}^{\tilde{r}_0}\tilde{c}_{\tilde{r}}(\varepsilon+n)^{\tilde{r}})\cdot$$
\begin{equation}\label{prodw}
\prod_{n=a-N}^{N-a}\prod_{i=1}^lw_i(\varepsilon-s_i+n)=
\end{equation}
$$=f_{1,N}(\varepsilon-N)\cdot f_{2,N}(\varepsilon+N)\cdot h_{N}(\varepsilon)\cdot\prod_{n=a-N}^{N-a}\prod_{i=1}^lw_{i,N,n}(\varepsilon-s_i+n)$$
where $f_1^0(\varepsilon),f_{1,N}(\varepsilon),f_2^0(\varepsilon),f_{2,N}(\varepsilon),w_i(\varepsilon),w_{i,N,n}(\varepsilon)$ are monic polynomials in $\varepsilon$. Notice that $h_1(\varepsilon),h_2(\varepsilon)$ can have poles in $\varepsilon$. 

We can simplify a bit the l.h.s. of (\ref{prodw}) using the following observation. The sum $\sum_{n=a-N}^{N-a}\sum_{\tilde{r}=0}^{\tilde{r}_0}\tilde{c}_{\tilde{r}}(\varepsilon+n)^{\tilde{r}}$ can be written as $p_1(\varepsilon-N)+p_2(\varepsilon+N)$ for some polynomials $p_1(\varepsilon),p_2(\varepsilon)\in {\bf k}[t]/(t^M)[\varepsilon]$ which are uniquely defined up to the transformations $p_1\mapsto p_1+c,~\p_2\mapsto p_2-c$ for $c\in {\bf k}[t]/(t^M)$.

The fact that the r.h.s. of (\ref{prodw}) has no poles in $\varepsilon$ implies that the l.h.s. can be written as 
$$f_1(\varepsilon-N)\cdot f_2(\varepsilon+N)\cdot\prod_{n=\tilde{a}-N}^{N-\tilde{a}}\prod_{i=0}^lw_i(\varepsilon-s_i+n)$$
for some $\tilde{a}\geqslant a$ independent of $N$, where $f_1,f_2$ are polynomials in $\varepsilon$. This implies that $w_{i,N,n}=w_i$ for $i=1,...,l$ and $\tilde{a}-N\leqslant n\leqslant N-\tilde{a}$, as well as factorization (\ref{detb}).

In this way, we obtain the statement of the proposition. 

$\square$

It follows from (\ref{rootsB}) and Proposition 3.1 that 
$$w_i(\lambda_{i,j})=0~~~\text{mod}~t^M$$
for $i=1,...,l$, $j=1,...,k_i$ and all sufficiently large $M$. Notice that the coefficients at $t^m,~m=1,...,M-1$ of $w_i(\varepsilon)$ are independent of $M$, so $w_i(\varepsilon)$ are defined as elements in ${\bf k}[[t]][\varepsilon]$. Therefore, 
$$w_i(\lambda_{i,j})=0$$
as a power series in $t$. This gives the following factorizations
$$w_i(\varepsilon)=\prod_{j=1}^{k_i}(\varepsilon-\lambda_{i,j}),~~~i=1,...,l.$$

Let us prove that the l.h.s. of (\ref{pol}) is equal to 
$$\prod_{i=1}^lw_i(\varepsilon)$$
which will give the statement of the theorem. 

Consider the subsum of the sum in the r.h.s. of (\ref{log}) with $n_1=0$. This subsum gives the matrix element $\log(1+A_N)_{0,0}$ which is equal to $\text{Coef}_{x^{\varepsilon}}(\log(1+P^{-1}Q)(x^{\varepsilon}))$ for sufficiently large $N$. Denote this matrix element by $R(\varepsilon)\in{\bf k}[\varepsilon,~\frac{1}{\varepsilon-s_i+m},~i=1,...,l,~m\in\Z]t{\bf k}[t]/(t^M)$.
Using partial fraction decomposition we can write
$$R(\varepsilon)=\sum_{r=1}^{r_0}\sum_{m=m_1}^{m_2}\frac{\bar{c}_{i,m,r}}{(\varepsilon-s_i+m)^r}+\sum_{s=0}^{s_0}\hat{c}_s\varepsilon^s$$
where $r_0,s_0\in\N~,m_1,m_2\in\Z$ are independent of $N$. We have 
\begin{equation}\label{log1}
\text{tr}~\log(1+A_N(\varepsilon)))=\sum_{n=-N}^NR(\varepsilon+n)+p_0(\varepsilon)
\end{equation}
where $p_0(\varepsilon)$ is a rational function in $\varepsilon$ which can have 
poles only at $\varepsilon=s_i-m$, $i=1,...,l$ and $|m-N|\leqslant b$ for some $b$ independent of $N$. On the other hand, comparing (\ref{log1}) and (\ref{det}) we get 
$$\sum_{m=m_1}^{m_2}\tilde{c}_{i,m,r}=c_{i,r}.$$
Therefore 
$$\varepsilon^ke^{T_{\varepsilon}(\text{Coef}_{x^{\varepsilon}}(\log(1+P^{-1}Q)(x^{\varepsilon})))}=\varepsilon^ke^{T_{\varepsilon}\Big(\sum_{r=1}^{r_0}\sum_{m=m_1}^{m_2}\frac{\bar{c}_{i,m,r}}{(\varepsilon-s_i+m)^r}+\sum_{s=0}^{s_0}\hat{c}_s\varepsilon^s\Big)}=$$
$$=\varepsilon^{k_1+...+k_l}e^{\sum_{r=1}^{r_0}\sum_{i=1}^l\sum_{m=m_1}^{m_2}\frac{c_{i,m,r}}{\varepsilon^r}}=\prod_{i=1}^l\varepsilon^{k_i}e^{\sum_{r=1}^{r_0}\frac{c_{i,r}}{\varepsilon^r}}=\prod_{i=1}^lw_i(\varepsilon).$$

This proves Theorem 3.2.

  $\square$

  {\bf Remark 3.2.} The fact that the l.h.s. of the formula (\ref{pol}) is a monic polynomial of degree $k$ has a shorter proof based on Theorem 2.2. It follows from the  definition (\ref{T}) of the operator $T_{\varepsilon}$ that 
  $$T_{\varepsilon}(\text{Coef}_{x^{\varepsilon}}(\log(1+P^{-1}Q)(x^{\varepsilon})))=\sum_{j=1}^{\infty}\frac{q_j(t)}{\varepsilon^j}~\in t\cdot {\bf k}[\varepsilon^{-1}][[t]]$$
  for some $q_j(t)\in {\bf k}[[t]]$ and the series in the r.h.s. is finite modulo any power of $t$. So the r.h.s. of (\ref{pol}) has a form
  $$\varepsilon^k+a_1\varepsilon^{k-1}+...+a_k+\frac{a_{k+1}}{\varepsilon}++\frac{a_{k+2}}{\varepsilon^2}+...$$
  for some $a_{l}(t)\in {\bf k}[[t]]$. On the other hand, if $r_1,...,r_k\in\Q$, then by Theorem 2.2 we have $a_j(t)=0~\text{mod}~p$ for $j>k$ and all large enough primes $p$. Indeed, $\det(\rho_{p,\xi,\eta}(D))$ is a polynomial in $\eta$. Therefore, $a_j(t)=0$ for $j>k$ and $r_1,...,r_k$ are rational. But $a_j(t)$ is a power series in $t$ and each coefficient is a rational expression in $r_1,...,r_k$. If this expression vanishes for infinitely many values of $r_1,...,r_k$, then it vanishes identically.

\section{Relation with monodromy}

In the case ${\bf k}\subset\C$ and the orders of all differential operators $Q_1,Q_2,...$ are not larger than $k$ we can write a formula for $w(\varepsilon)$ in terms of a certain monodromy operator.

Fix $M\geqslant 1$ and consider a cyclic module $\mathcal{E}$ over the ring $R={\bf k}(\!(x)\!)[\partial_x][t]/(t^M)$ generated by an element $\psi$ such that $(P+Q)\psi=0$, i.e. $\mathcal{E}\cong R/R (P+Q)$. The assumption that the order of $P+Q$ is $k$ implies that $\mathcal{E}$ is a free ${\bf k}(\!(x)\!)[t]/(t^M)$-module of rank $k$. Hence, $\mathcal{E}$ has rank $kM$ as a ${\bf k}(\!(x)\!)$-module.
 We can interpret it as the space of sections of a vector bundle $\mathcal{V}$ 
 over the formal punctured disc. The structure of the ${\bf k}(\!(x)\!)[\partial_x]$-module gives a meromorphic flat connection $\nabla$ on $\mathcal{V}$. Notice that $\mathcal{V}$ can also be considered as a bundle of free ${\bf k}[t]/(t^M)$-modules of rank $k$.

 This connection has regular singularities, that is, after a gauge transformation from $GL_{kM}({\bf k}(\!(x)\!)$, it will have only simple poles at $x=0$. Indeed, as a 
 ${\bf k}(\!(x)\!)[\partial_x]$-module $\mathcal{E}$ has finite filtration $\mathcal{E}\supset t\mathcal{E}\supset...\supset t^{M-1}\mathcal{E}$ with associated graded modules given by equation $P\psi=0$. According to the Hukuhara-Levelt-Turrittin formal classification of meromorphic connections \cite{A}, any extension of formal meromorphic connections with regular singularities is also regular. 

 Finally, after the extension of scalars from ${\bf k}$ to $\C$, we can identify our formal meromorphic connection with regular singularities with analytic meromorphic connection on a punctured disc $\{x\in\C;~0<x<r\}$ for some $r>0$. For this analytic connection, we have a monodromy operator $\mathcal{M}\in
  GL_k(\C[t]/(t^M)$. The eigenvalues of $\mathcal{M}$ are $e^{2\pi i(s_i+\lambda_{i,j})}$, $s_i\in\C$, $\lambda_{i,j}\in t\C[t]/(t^M)$.

  Let $\log \mathcal{M}$ be the value of $\log$ of the operator $\mathcal M$ with the eigenvalues $\{2\pi i(s_i+\lambda_{i,j})\}$. It commutes with the operator $\mathcal{S}$ which, for each $i=1,...,l$, acts by  multiplication by $s_i$ on the space of generalized eigenvectors of the operator 
  $\log\mathcal{M}$ with eigenvalue $\{2\pi i(s_i+\lambda_{i,j})\}$. Then we have \begin{equation}\label{mon}
  w(\varepsilon)=(-1)^k\det \Big (\frac{1}{2\pi i}\log \mathcal{M}-\mathcal{S}-\varepsilon\Big).
  \end{equation}
  Moreover, if $P+Q$ is convergent for $t,x\in\C$ such that $|t|<r$, $0\leqslant |x|\leqslant r$ for some $r>0$, then $\mathcal{M}$ will be an analytic function in $t$, whose Taylor expansion up to $t^{M-1}$ coincides with (\ref{mon}).

  {\bf Remark 4.1.} For a differential operator $D$ with symbol $\partial_{\theta}^k$ on a circle $S_{\theta}^1$, where $\theta=\text{arg}~x$, one can define a $\zeta$-regularized determinant. There exists a closed formula \cite{BF} 
  $${\det}_{\zeta}D=\exp(\phi(D))\cdot\det(1-\mathcal{M})$$
  where $\mathcal{M}$ is the monodromy of the local system of solutions $D\psi=0$, and $\phi(D)$ is an explicit functional of $D$. The operator $D$ also acts on sections of the line bundle with flat connection and monodromy $e^{2\pi i\varepsilon}$ for any $\varepsilon\in\C$. Let us denote by $D_{\varepsilon}$ the resulting operator. The corresponding $\zeta$-regularized determinant ${\det}_{\zeta}D$ up to an invertible element is 
  equal to 
  $$\det(\mathcal{M}-e^{2\pi i\varepsilon}).$$
  This formula is curiously similar, but not equivalent to our formula (\ref{mon}).

  \section{The case of \texorpdfstring{$q$}{q}-difference operators}

  Let $A={\bf k}\langle q^{\pm 1},x^{\pm 1},y^{\pm 1}\rangle$ be an algebra\footnote{This algebra is related with algebra of differential operators by $y=\exp(\hbar x\partial_x),~q=\exp(\hbar).$} generated by $q,x,y$ with defining relations $yx=qxy$ and $q\in A$ is a central element.

  Let $\bar{A}={\bf k}\langle q^{\pm 1},x^{\pm 1},y^{\pm 1},(y-q^j)^{-1},~j\in\Z\rangle$. An arbitrary element $S\in\bar{A}$ can be written uniquely in the form
  $$S=\sum_{j=m_1}^{m_2}x^jf_j(y)$$
  where $m_1,m_2\in\Z$ and $f_j(y)\in{\bf k}[q^{\pm 1}](y)$ are rational functions in $y$, such that denominators of all $f_j(y)$ are monic polynomials with roots belonging to the set $\{q^j;~~j\in\Z\}$. 
  We can multiply such expressions using commutation relations
  $$f(y)x^j=x^jf(q^jy)$$
  for an arbitrary rational function $f(y)$. 

  For each natural number $n\geqslant 1$ we define the algebra $A_n$ as a quotient
  $$\bar{A}_n=\bar{A}/\Phi_n(q)$$
  where $\Phi_n(q)\in\Z[q]$ stands for $n$-th cyclotomic polynomial.

  For each $n\geqslant 1$ and formal parameters $\xi,\eta$ we define an $n$-dimensional matrix representation $\rho_{n,\xi,\eta}$ of the algebra $\bar{A}$, which is also a representation of its quotient $\bar{A}_n$, 
  as follows
\begin{equation}\label{rhon}
	\rho_{n,\xi,\eta}(x)=	\begin{pmatrix}
			0 & 0 & & \xi \\
			\xi & 0 & \ddots & \\
			& \ddots & \ddots & 0 \\
			0 & & \xi & 0
		\end{pmatrix}, \quad
		\rho_{n,\xi,\eta}(y)=\begin{pmatrix}
			\eta & 0 & & 0 \\
			0 & q\eta & \ddots & \\
			& \ddots & \ddots & 0 \\
			0 & & 0 & q^{n-1}\eta
		\end{pmatrix}.
	\end{equation}
    The corresponding action in the basis $\{v_j;~j\in\Z/(n)\}$ is given by 
    $$\rho_{n,\xi,\eta}(x)v_j=\xi v_{j+1},~~~\rho_{n,\xi,\eta}(y)v_j=q^j\eta v_j.$$

   For a formal variable $q^{\varepsilon}$ we define\footnote{In the sequel we will use the notations $q^{\varepsilon+m}=q^{\varepsilon}q^m$, ~~~$q^{m\varepsilon}=(q^{\varepsilon})^m$ for arbitrary $m\in\Z$.} an infinite-dimensional representation $\rho_{\varepsilon}$ of the algebra $\bar{A}$ in the vector space over ${\bf k}(q,q^{\varepsilon})$ with the basis 
   $\{x^{\varepsilon+m};~m\in\Z\}$ by
   $$\rho_{\varepsilon}(x):~~~x^{\varepsilon+m}\mapsto x^{\varepsilon+m+1},$$
   $$\rho_{\varepsilon}(y):~~~x^{\varepsilon+m}\mapsto q^{\varepsilon+m}x^{\varepsilon+m}.$$
   In other words, for arbitrary element $f(x)\in x^{\varepsilon}{\bf k}(q,q^{\varepsilon})[x^{\pm 1}]$ we define
   \begin{equation}\label{qact}
   \rho_{\varepsilon}(x):~~~f(x)\mapsto xf(x),~~~~~~~~~\rho_{\varepsilon}(y):~~~f(x)\mapsto f(qx).
   \end{equation}

  Define a linear operator 
  $$R_{n,\varepsilon}:~{\bf k}[q^{\pm 1},q^{\pm\varepsilon},(1-q^{\varepsilon+j})^{-1},~j\in\Z]\to {\bf k}[q^{\pm 1},q^{n\varepsilon},(1-q^{n\varepsilon})^{-1},~ (q^i-1)^{-1},~i\geqslant 1]$$
  by
  $$R_{n,\varepsilon}(q^{j\varepsilon})=0,~j\ne 0,~~~~R_{n,\varepsilon}(1)=n,~~~~~~R_{n,\varepsilon}\Bigg(\frac{1}{(q^{\varepsilon+i}-1)^j}\Bigg)=\frac{nW_j(n,q^{n\varepsilon})}{(q^{n\varepsilon}-1)^j},~j > 0,~i\in\Z$$
  where we use partial fractions decomposition\footnote{Notice that polynomials $q^i-1$ appear in denominators after this decomposition.} with respect to variable $q^{\varepsilon}$, and $W_j(n,q^{n\varepsilon})$ is a  polynomial in $q^{n\varepsilon}$ of degree $j-1$, its coefficients are polynomials in $n$ having integer values for integer $n$. These polynomials are
  defined by 
  $$\frac{nW_j(n,q^{n\varepsilon})}{(q^{n\varepsilon}-1)^j}=\sum_{i=0}^{n-1}\frac{1}{(\mu^iq^{\varepsilon}-1)^j}$$
  where $\mu$ is a prime $n$-th root of unity. We also have an equation in $\Q[n][[q^{n\varepsilon}]]$
  $$W_j(n,q^{n\varepsilon})=(1-q^{n\varepsilon})^j\sum_{i=0}^{\infty}\frac{(ni+1)...(ni+j-1)}{(j-1)!}q^{ni\varepsilon}.$$
  For example
  $$W_1(n,u)=1,~~~W_2(n,u)=1+(n-1)u,~~~W_3(n,u)=1+\frac{1}{2}(n-1)(n+4)u+\frac{1}{2}(n-1)(n-2)u^2,$$
  $$W_4(n,u)=1+\frac{1}{6}(n-1)(n^2+7n+18)u+\frac{1}{3}(n-1)(2n^2+2n-9)u^2+\frac{1}{6}(n-1)(n-2)(n-3)u^3.$$

  In particular, let $R_{\varepsilon}=R_{1,\varepsilon}$, then
  $$R_{\varepsilon}(q^{j\varepsilon})=0,~j\ne 0,~~~~R_{\varepsilon}(1)=1,~~~~~~R_{\varepsilon}\Bigg(\frac{1}{(q^{\varepsilon+i}-1)^j}\Bigg)=\frac{1}{(q^{\varepsilon}-1)^j},~j > 0,~i\in\Z.$$
  Notice that for an arbitrary $f\in {\bf k}[q^{\pm 1},q^{\pm\varepsilon},(1-q^{\varepsilon+j})^{-1},~j\in\Z]$ we have
  \begin{equation}\label{rnr1}
  R_{n,\varepsilon}(f)=\sum_{j=0}^{n-1}R_{\varepsilon}(f)|_{q^{\varepsilon}~\to ~q^{\varepsilon+j}}~~~\text{mod}~\Phi_n(q).
  \end{equation}

  Introduce the notation
  \begin{equation}\label{Pq}P=\prod_{j=1}^k(y-q^{l_j}),~~~l_j\in \Z
  \end{equation}

  {\bf Theorem 5.1.} Let $U=P^{-1}D_1P^{-1}D_2...D_NP^{-1}\in\bar{A}$ 
  where $D_1,...,D_N\in A$. Then for all large enough numbers $n>C$ where $C$ depends only on $U$  we have
  $$\text{tr}~ (\rho_{n,\xi,\eta}(U))=R_{n,\varepsilon} \text{Coef}_{x^{\varepsilon}}(\rho_{\varepsilon}(U)x^{\varepsilon})|_{q^{\varepsilon}=\eta} ~~~\text{mod}~\Phi_n(q).$$

  {\bf Theorem 5.2.}  Let $D=P+Q$ where $P$ is given by (\ref{Pq}) and $Q\in t{\bf k}[q^{\pm 1}]((x))[y^{\pm 1}][[t]]$. Then for any $M\geqslant 1$ there exists $C_M$ such that
  $$\det(\rho_{n,\xi,\eta}(D))=(-1)^{k(n-1)}(q^{n\varepsilon}-1)^k\cdot e^{R_{n,\varepsilon}\text{Coef}_{x^{\varepsilon}}(\rho(\log(1+P^{-1}Q))x^{\varepsilon})}|_{q^{\varepsilon}=\eta}~~~\text{mod}~\Phi_n(q),~t^{M}$$
  for all  $n>C_M$.

  {\bf Theorem 5.3.} Define the action of the algebra $\bar{A}[[t]]$ on the space 
  ${\bf k}(q)(\!(x)\!)[\log_q x][[t]]$ by the formulas (\ref{qact}). In particular, 
  $y\log_q x=\log_q x+1$. Then the space of solutions of the equation
  $$(P+Q)f=0,~~~f\in{\bf k}(\!(x)\!)[\log_q x][[t]]$$
  is a free $k$-dimensional module over ${\bf k}[[t]]$. Moreover,  after an  extension of scalars the module of solutions has a basis 
  $$f_1,...,f_k\in x^{\lambda_i}\overline{{\bf k}(q)}(\!(x)\!)[[t^{1/L}]][\log_q x]$$
  for some $L\geqslant 1$, where $q^{\lambda_i}-1\in t^{1/L}\overline{{\bf k}(q)}[[t^{1/L}]]$ for $i=1,...,k$ so we can write $\lambda_i\log q\in t^{1/L}\overline{{\bf k}(q)}[[t^{1/L}]]$, degrees of $f_i$ as polynomials in $\log_q x$ are smaller than $k$, and $x^{\lambda_i}$ are understood as
  $$x^{\lambda_i}=e^{\lambda_i\log q \log_q x}=\sum_{m=0}^{\infty}\frac{(\lambda_i\log q)^m(\log_q x)^m}{m!}~\in \overline{{\bf k}(q)}(\!(x)\!)[\log_q x][[t^{1/L}]].$$

  {\bf Theorem 5.4.} There exist a monic polynomial $$w(q^{\varepsilon})=q^{k\varepsilon}+u_1q^{(k-1)\varepsilon}+...+u_k\in {\bf k}[q^{\pm 1},(q^i-1)^{-1},~i\geqslant 1][[t]][q^{\varepsilon}]$$
  such that $w(q^{\varepsilon})=(q^{\varepsilon}-1)^k$ mod $t$ and,  we have 
  \begin{equation}\label{qw}
  w(q^{\varepsilon})=\prod_{i=1}^k(q^{\varepsilon}-q^{\lambda_i})
  \end{equation}
  where $q^{\lambda_1},...,q^{\lambda_k}$ are defined in Theorem 5.3. Then be have
\begin{equation}\label{R1w}(q^{\varepsilon}-1)^k\cdot e^{R_{\varepsilon}\text{Coef}_{x^{\varepsilon}}(\rho(\log(1+P^{-1}Q))x^{\varepsilon})}=C w(q^{\varepsilon})
\end{equation}
where $C\in \in {\bf k}[q^{\pm 1},(q^i-1)^{-1},~i\geqslant 1][[t]]$ and $C=1$ mod $t$.
Moreover, if degree of $Q$ with respect to $y$ is smallet than $k$, then $C=1$. In this case we also have 
\begin{equation}\label{Rnw}
(q^{n\varepsilon}-1)^k\cdot e^{R_{n,\varepsilon}\text{Coef}_{x^{\varepsilon}}(\rho(\log(1+P^{-1}Q))x^{\varepsilon})}=\prod_{i=0}^{k}(q^{n\varepsilon}-q^{n\lambda_i})~~~\text{mod}~\Phi_n(q),~t^M
\end{equation}
for $n>C_M$.

If all coefficients of $P+Q$ belong to a subring $R\subset {\bf K}[q^{\pm 1}][[t]]$, then $u_1,...,u_k\in R[q^{\pm 1},~(q^i-1)^{-1},~i\geqslant 1][[t]]$.

\

  {\bf Proof.} Proofs of Theorems 5.1, 5.2, 5.3 are similar to ones of Theorems 2.1, 2.2, 3.1.
   Let us outline the proof of Theorem 5.4.

  Similarly to the proof of Theorem 3.2 we assume without loss of generality that $L=1$, the field ${\bf k}$ is algebraically closed, and all $\lambda_1,...,\lambda_k$ are distinct. Fix sufficiently large $M>1$ and work modulo $t^M$. 
 Let $B(q^{\varepsilon})$ be an infinite matrix corresponding to the action of $P+Q$  on the free ${\bf k}(q,q^{\varepsilon})[[t]]$-module with the basis $\{x^{n+\varepsilon};~n\in\Z\}$.  For $N\geqslant 1$ we denote by $B_{N}(q^{\varepsilon})$ the finite submatrix of $B(q^{\varepsilon})$ corresponding to indices $n\in \{-N,...,N\}$.

 We have $B_N(q^{\lambda_i})f_i=0$ for $i=1,...,k$ and therefore
 $$\det B_N(q^{\lambda_i})=0,~~~i=1,...,k.$$

 {\bf Proposition 5.1.} There exist $a\geqslant 0$ depending only on $P+Q$, such that if $N\geqslant a$, then we have 
 \begin{equation}\label{detb1}\det B_{N}(q^{\varepsilon})=(1+c_1(q^{\varepsilon-N}))(1+c_2(q^{\varepsilon+N}))f_1(q^{\varepsilon-N})~\prod_{j=a-N}^{N-a}w(q^{\varepsilon+j})~f_2(q^{\varepsilon+N}).
 \end{equation}
 Here $f_1(q^{\varepsilon}),f_2(q^{\varepsilon}),w(q^{\varepsilon}),c_1(q^{\varepsilon}),c_2(q^{\varepsilon})\in R[\frac{1}{q^i-1},i=1,2,...][q^{\varepsilon}]$ are such that

 \begin{align}1)&~~~ f_1(q^{\varepsilon})=f_1^{(0)}(q^{\varepsilon})\text{~~~mod~}t\\
2)&~~~ f_2(q^{\varepsilon})=f_2^{(0)}(q^{\varepsilon})\text{~~~mod~}t\\
 3)&~~~w(q^{\varepsilon})=(q^{\varepsilon}-1)^k\text{~~~mod~}t\\
 4)&~~~c_1(q^{\varepsilon}),~c_2(q^{\varepsilon})=0\text{~~~mod~}t\end{align}
where 
$$f_1^{(0)}(q^{\varepsilon})=\prod_{\substack{1\leqslant i\leqslant k\\0\leqslant j<l_i+a}}(q^{\varepsilon-l_i+j}-1),~~~~~~~f_2^{(0)}(\varepsilon)=\prod_{\substack{1\leqslant i\leqslant k\\l_i+a< j<0}}(q^{\varepsilon-l_i+j}-1)$$
and $f_1(q^{\varepsilon}),f_2(q^{\varepsilon}),w(q^{\varepsilon})$ are monic polynomials.

We use here the same Lemma 3.1, and the analog of Lemma 3.2 is

{\bf Lemma 5.1.}  Let $A(q^{\varepsilon})=(A_{i,j}(q^{\varepsilon}))_{i,j\in\Z}$ be an infinite matrix such that its matrix elements have a form $A_{i,j}(q^{\varepsilon})=a_{j-i}(q^{\varepsilon+i})$ where $a_i(q^{\varepsilon})\in {\bf k}[q^{\pm 1},q^{\varepsilon},\frac{1}{q^{\varepsilon+\Z}-1}][t]/(t^M)$ and $a_i(q^{\varepsilon})=0$ for $|i|\gg 0$. Then for submatrix of $1+A(q^{\varepsilon})$ corresponding to indices $n\in\{-N,...,N\}$ we have 
$\log\det(1+A_N(q^{\varepsilon}))\in {\bf k}[q^{\varepsilon},\frac{1}{q^{\varepsilon+\Z}-1}][t]/(t^M)$ and for $N\gg 0$
\begin{equation}\label{qdet}
\log\det(1+A_N(q^{\varepsilon}))=h_1(q^{\varepsilon-N})+\sum_{n=a-N}^{N-a}\Big(\sum_{r=1}^{r_0}\frac{c_r}{(q^{\varepsilon+n}-1)^r}+\sum_{\tilde{r}=0}^{\tilde{r}_0}\tilde{c}_{\tilde{r}}(q^{\varepsilon+n}-1)^{\tilde{r}}\Big)+h_2(q^{\varepsilon+N})
\end{equation}
for some $a,r_0,\tilde{r}_0\geqslant 0$, where $c_r,\tilde{c}_r\in t{\bf k}[t]/(t^M)$ and $h_1(\varepsilon),h_2(\varepsilon)\in t{\bf k}[q^{\varepsilon},\frac{1}{q^{\varepsilon+\Z}-1}][t]/(t^M)$. Using these lemmas, we prove the proposition.

$\square$

It follows from the Proposition 5.1 that $w(q^{\varepsilon})$ is defined as an element in ${\bf k}[q^{\pm 1},(q^i-1)^{-1},~i\geqslant 1][[t]][q^{\varepsilon}]$ and given by (\ref{qw}). The formula (\ref{R1w}) can be proved in the same way as
(\ref{pol}). To prove (\ref{Rnw}) notice that it follows from (\ref{rnr1}) that
$$\prod_{i=1}^{n-1}w(q^{\varepsilon+j})=\prod_{j=0}^{n-1}(q^{\varepsilon+j}-1)^k\cdot e^{R_{n,\varepsilon}\text{Coef}_{x^{\varepsilon}}(\rho(\log(1+P^{-1}Q))x^{\varepsilon})}~~~\text{mod}~\Phi_n(q),~t^M$$
for $n>C_M$ and we can factorize $w$ using (\ref{qw}) and apply to both l.h.s. and r.h.s. of this identity the following formula
$$\prod_{j=0}^{n-1}(q^{\varepsilon+j}-\lambda)=(-1)^{n-1}(q^{n\varepsilon}-\lambda^n)~~~\text{mod}~\Phi_n(q).$$

Finally, the last assertion of the theorem concerning the coefficients $u_i$ of the polynomial $w(q^{\varepsilon})$, follows from the Remark 3.1.

$\square$

{\bf Remark 5.1.} The fact that the l.h.s. of formula (\ref{Rnw}) is a polynomial in $q^{n\varepsilon}$ has a shorter proof based on the result of Theorem 5.2, and similar to the Remark 3.2. It follows from definition of $R_{n,\varepsilon}$ that 
$$(q^{n\varepsilon}-1)^k\cdot e^{R_{n,\varepsilon}\text{Coef}_{x^{\varepsilon}}(\rho(\log(1+P^{-1}Q))x^{\varepsilon})}=\sum_{j=0}^{\infty}\psi_j~(q^{n\varepsilon}-1)^{k-j}$$
where $\psi_j$ is a power series in $t$ for each $j$. Each coefficient of this
power series of $t$ has the form $\frac{\mu_j(n,q)}{\Phi_{m_1}(q)...\Phi_{m_l}(q)}$
where $\mu_j(n,q)$ are polynomials in $n,q$. From Theorem 5.2 it follows that 
$\mu_j(n,q)=0$ mod $\Phi_n(q)$ for $j>k$ and $n\gg 1$, which means that
$\mu_j(n,q)=0$ for $j>k$.

{\bf Remark 5.2.} Using Theorem 5.4 we can reformulate the result of Theorem 5.2 (in the case of degree of $Q$ with respect to $y$ is less than $k$) as follows
$$\det(\rho_{n,\xi,\eta}(D))=(-1)^{k(n-1)}\prod_{i=0}^{k}(q^{n\varepsilon}-q^{n\lambda_i})|_{q^{\varepsilon}=\eta}~~~\text{mod}~\Phi_n(q),~t^M$$
for $n>C_M$.

\section{Explicit formulas in a certain class of deformations of identity operator}

In this section, we assume $A={\bf K}\langle q^{\pm1},x^{\pm 1},y^{\pm 1}\rangle$ with the same commutation relations as before. Let us consider the simplest nonzero case $P=1$. The results of Theorem 5.2 can be applied for the computation of $\det \rho_{n,\xi,\eta}(1+tQ)$ mod $t^M,~\Phi_n(q)$ for $n>C_M$.  However, there exists a different approach that gives an exact analytic formula for determinants without taking a truncation modulo $t^M$ for a certain class of $Q\in A$. It will be convenient to set $t=\lambda^{-1}$ and multiply our difference operator by
$\lambda$. We obtain the following family of difference operators
$$\lambda+Q$$
where $\lambda$ is a parameter. In the sequel, we will work in formal power series in $\lambda^{-1}$, so informally we assume $\lambda\gg 0$. We want to compute $\det(\rho_{n,\xi,\eta}(\lambda+Q))$ mod $\Phi_n(q)$.

Let us start with a particular example, and after that we will briefly outline the general result.

Consider a family  of $q$-difference operators
$$D=\lambda+(-x-y+x^{-1}y^{-1})\in A.$$

{\bf Theorem 6.1.} We have
\begin{equation}\label{D1}
\det(\rho_{n,\xi,\eta}(D))=-\xi^n-\eta^n+\xi^{-n}\eta^{-n}+F(G_q(\lambda)^n)~~~\text{mod}~\Phi_n(q)
\end{equation}
where 
$$F(u)=u-2u^{-2}+5u^{-5}-32u^{-8}+286u^{-11}-...$$
$$G_q(\lambda)=\lambda+(q^{-1}+1)\lambda^{-2}-(q^{-3}+3q^{-2}+5q^{-1}+3+q)\lambda^{-5}+...$$
are infinite series defined as follows:
$$G_q(\lambda)=\lambda\exp(\text{Coeff}_{0,0}(\log(\lambda^{-1}D))),
$$
$$F(G_1(\lambda))=\lambda.$$
Here $\exp$ and $\log$ are defined as infinite power series and $\text{Coeff}_{i,j}$ stands for a coefficient at $x^iy^j$ of a Laurent series in $x,~y$.

Notice that 
$$G_1(\lambda)=\lambda\exp\Big(\sum_{m=1}^{\infty}(-1)^{m-1}\frac{(3m-1)!}{(m!)^3}\lambda^{-3m}\Big).$$

{\bf Proof.} It is known that
$$\det(\rho_{n,\xi,\eta}(D))=-\xi^n-\eta^n+\xi^{-n}\eta^{-n}+H_q(\lambda)~~~\text{mod}~~\Phi_n(q)$$
where $H_q(\lambda)$ is a polynomial in $\lambda$ of degree $n$ with coefficients in $\Z[q]/\Phi_n(q)$. We have
$$\log(-\xi^n-\eta^n+\xi^{-n}\eta^{-n}+H_q(\lambda))=\log(\det(\rho_{n,\xi,\eta}(D))=$$
\begin{equation} \label{det1}
=\text{tr}~(\log(\rho_{n,\xi,\eta}(D))=\text{tr}~(\rho_{n,\xi,\eta}(\log(D))).
\end{equation}
The l.h.s. of (\ref{det1}) can be written as 
$$\log(H_q(\lambda))+\log(1+H_q(\lambda)^{-1}(-\xi^n-\eta^n+\xi^{-n}\eta^{-n}))$$
and coefficient of this series at $x^0y^0$ is equal to 
$$\log(H_q(\lambda))+\text{Coeff}_{0,0}(\log(1+H_q(\lambda)^{-1}(-\xi^n-\eta^n+\xi^{-n}\eta^{-n})))=\log(G_1(H_q(\lambda)))$$
The r.h.s. of (\ref{det1}) can be written as
$$\text{tr}~(\rho_{n,\xi,\eta}(\log(\lambda)+\log(\lambda^{-1}D)))$$
and coefficient of this series at $x^0y^0$ is equal to
$$n\log(\lambda)+n\text{Coeff}_{0,0}(\log(\lambda^{-1}D))=\log(G_q(\lambda)^n).$$
Equating exponentials we get the identity
$$G_1(H_q(\lambda))=G_q(\lambda)^n~~~\text{mod}~~\Phi_n(q)$$
which is equivalent to the statement of the theorem. 

$\square$

{\bf Remark 6.1.} This result can be generalized to the case of difference operators of the form
$$D=\sum_{(i,j)\in\Delta}c_{i,j}x^iy^j$$
where $\Delta\in\R^2$ is a convex polygon with integer vertices and one integer interior point. Without loss of generality we can assume that the interior point of $\Delta$ is $(0,0)\in\Delta$ and write our difference 
operator as 
\begin{equation}\label{Dg}
D=\lambda+\sum_{(i,j)\in\bar{\Delta}}c_{i,j}x^iy^j
\end{equation}
where $\bar{\Delta}$ is the boundary of $\Delta$. In this case, the formula
for $\det(\rho_{n,\xi,\eta}(D))$ and its proof are essentially the same.

{\bf Example 6.1.} Let 
$$D=\lambda-x-y-ax^{-1}-by^{-1}$$
where $a,b$ are parameters. Then we have
$$\det(\rho_{n,\xi,\eta}(D))=-\xi^n-\eta^n-a^n\xi^{-n}-b^n\eta^{-n}+F(a^n,b^n,G_q(a,b,\lambda)^n)~~~\text{mod}~\Phi_n(q)$$
where $G_q(a,b,\lambda),~F(a,b,u)$ are defined by
$$G_q(a,b,\lambda)=\lambda\exp(\text{Coeff}_{0,0}(\log(\lambda^{-1}D))),
$$
$$F(a,b,G_1(a,b,\lambda))=\lambda.$$

{\bf Remark 6.2.} In general, if $(0,0)\in\Delta$ is the only interior integer point, the formulas for $\det(\rho_{n,\xi,\eta}(D))$ are similar, we just need to choose a correct parameterization of the boundary part $\sum_{(i,j)\in\bar{\Delta}}c_{i,j}x^iy^j$ of our difference operator $D$. 
Namely, for each side of the polygon $\Delta$ the corresponding subsum of
(\ref{Dg}) has a form 
$$\sum_{r=0}^{s}c_{i_0+ri_1,j_0+rj_1}x^{i_0+ri_1}y^{j_0+rj_1}$$
where $i_1,j_1$ are coprime. We need to write this in a factorized form
$$c_{i_0,j_0}x^{i_0}y^{j_0}(1+\mu_1x^{i_1}y^{j_1})...(1+\mu_sx^{i_1}y^{j_1})$$
and take $(-1)^{i_0+j_0}c_{i_0,j_0},(-1)^{i_1+j_1}\mu_1,...,(-1)^{i_1+j_1}\mu_s$ for all sides of $\bar{\Delta}$ as new parameters.

\section{Adic properties of polynomial \texorpdfstring{$w$}{w} in some examples}

In this section, we reproduce some experimental results for a family of differential operators discussed also in \cite{KO}, and its $q$-analog. We also describe some curious adic properties of the polynomial $w$ in these
examples. It turns out that the hypothetical formula for factorization of the denominators of the coefficients of $t$-expansion is much simpler for $q$-analog, and a more complicated formula in the arithmetic case can be obtained from the formula for $q$-analog.

\subsection{Arithmetic case}

Consider the differential operator depending on the parameter $t$:

$$D_t=\partial_xx(x-1)(x-t)\partial_x+x$$
It has regular singularities at $0,1,t,\infty$.  

It acts in the basis $(x^n)_{n\in \Z}$ as
\begin{equation}\label{mdt}
\begin{array}{c|cc}\, & x^{n-1}&x^n\\ \hline x^{n-1}& (n-1)n(t+1) & n^2 t\\
x^n & n^2 & n(n+1)(t+1)\end{array}
\end{equation}
In fact, the operator $D_t$ at $t=0$ does not belong to the family (\ref{P}).
However, in the rescaled basis $(t^{n/2}x^n)_{n\in \Z}$ it acts as a symmetric matrix
$$\begin{array}{c|cc}\, & t^{(n-1)/2}x^{n-1}&t^{n/2}x^n\\ \hline  t^{(n-1)/2}x^{n-1}& (n-1)n(t+1) & n^2 t^{1/2}\\
t^{n/2}x^n & n^2t^{1/2} & n(n+1)(t+1)\end{array}$$
i.e. operator became
$$ x\partial_x(x\partial_x+1)  +t^{1/2}\big[x^{-1}(x\partial_x)^2+(x\partial_x)^2x\big]+tx\partial_x(x\partial_x+1) $$
which is a deformation of $x\partial_x(x\partial_x+1).$

For small $t\in \C$, consider the monodromy operator $\mathcal M_t$ of the local system of solutions of $D_t$ around a loop surrounding singular points $0,t$. The eigenvalues of $\mathcal M_t$ are of the form $\exp(\pm 2\pi i \lambda(t))$ where $\lambda(t)={1\over 2 }t+\dots \in \Q[[t]]$ is a convergent series.
The characteristic polynomial of ${1\over 2\pi i}\log(\mathcal M_t)$ is 
$$w(\varepsilon)= \det \Big(\varepsilon -{1\over 2\pi i}\log(\mathcal M_t)\Big )=\varepsilon^2-\lambda(t)^2 $$
Let us consider series
$$-\lambda(t)^2=c_1 t+c_2 t^2+\dots=-{1\over 4} t^2-{1\over 24} t^3-{101\over 576 } t^4+\dots$$
We calculated coefficients $c_n$ for $n\leqslant 251$.

Using these data, we found:

$$c_n=\frac{b_n}{\prod_{p\leqslant n}p^{\alpha_p(n)}},~~~n\geqslant 2$$
\begin{equation}\label{cn}
\alpha_2(n)=\text{ord}_2(2^{n-1}n!),~~~2\leqslant n,
\end{equation}
$$\alpha_p(n)=\text{max}_{l\geqslant 1}~l(n+1-p^l),~~~3\leqslant p\leqslant n$$
$$\text{with~exceptions~~~}\alpha_3(6)=2,~\alpha_3(7)=4,~\alpha_3(8)=5$$
where $p$ stands for prime number, and numerator and denominator in the formular for $c_n$ are coprime integers.
Moreover, the following properties hold.

{\bf Claim 7.1} For $p<n<2 p^2-p$  where prime $p$ lies between 13 and 250 one has	$${c_{n+1}\over c_n}= -p^{-1}\cdot g(p)$$
$$g(p)=\Big(1-{1\over 12}p+{49\over 144}p^2+{751\over 8640}p^3+{28777\over518400}p^4+{6903793\over 21772800}p^5+\dots  \Big) $$

Denominators of coefficients of $g$:

$2^2 ~ 3^1$

$2^4~ 3^2 $

$2^6 ~3^3~5^1$

$2^8 ~3^4 ~ 5^2$

$2^{10} ~ 3^5 ~5^3 ~ 7^1$

If we go up to $p^6$: new phenomenon, the ratio $c_{n+1}/ c_n$ does not stabilize, depends on $n\mod 2$.

Add to $g(p)$ term
$${1686498937\over 9144576000} p^6\quad \implies \quad {c_{n+2}\over c_n}={1\over p^2}\cdot g(p)^2\cdot (1+O(p^7))$$
Denominator of the last coefficient is

 $ 2^{12}~ 3^6~ 5^4 ~ 7^2$

up to $p^7$: two different truncated series

$$g_{even}(p)=
1-
{1\over 12}p+{49\over 144}p^2+{751\over 8640}p^3+{28777\over 518400}p^4+{6903793\over 217728000}p^5+{1686498937\over 91445760000}p^6+{210013626433\over 38407219200000}p^7$$

$$g_{odd}(p)=
1-
{1\over 12}p+{49\over 144}p^2+{751\over 8640}p^3+{28777\over 518400}p^4+{6903793\over 217728000}p^5+{1686498937\over 91445760000}p^6+{160004226433\over 38407219200000}p^7$$

(last coefficients differ by $1/768$).

Then
 $${c_{n+2}\over c_n}={1\over p^2}\cdot g_{even}(p)^2\text { or }{1\over p^2}\cdot g_{odd}(p)^2 \quad \text{    depending on }\quad n\mod 2$$

The denominator of the coefficient of $p^7$ in both series is

 $2^{14}~ 3^7 ~ 5^5 ~7^3$
 
General formula for the denominator for the coefficient of $p^k$:
 $$  2^{2k} \cdot 3^{3+k-3}\cdot 5^{3+k-5}\cdot 7^{3+k-7}\cdot\dots$$

{\bf Corollary 7.1.} $p$-adic number 

$$c_n\cdot (-1)^nn \cdot p^{n+1-p}\cdot g(p)^{-n}$$
does not depend on $n$ up to $p^6$:
It is equal to
 $$ 1-{3\over 2}p-{19\over 24}p^2 +O(p^3). $$

 {\bf Remark 7.1.} For a given prime $p\geqslant 3$, the $p$-adic properties of $c_n$ 
 described by (\ref{cn}) suggest that analytic properties of the series 
 $\sum_{n\geqslant 2}c_nt^n$ as $p$-adic function in $t$ can be modeled by the anzats
 \begin{equation}\label{apr}
 \sum_{n\geqslant 2}c_nt^n~\sim~\sum_{l=1}^{\infty}\frac{a_l}{t-b_l},~~~a_l,b_l\in\Q_p,~~~|a_l|_p=p^{-lp^l},~|b_l|_p=p^{-l}.
 \end{equation}
 Notice that any such series gives a well-defined function on $\C_p\setminus \{0,b_m,~m=1,2,...\}$, therefore, the series $\sum_{n\geqslant 2}c_n t^n$ is divergent. The experimental data suggest that (\ref{apr}) approximates well the behavior of the coefficients $c_n$ in the interval $p<n<2p^2-p$. However, looking with sufficiently high precision, we see a deviation from the geometric series. 
 Presumably, we have a well-defined analytic function in a neighbor of 0 with removed $\{0\}$ and a disjoint union of small discs instead of punctures. The approximate position $b_1$ of the first pole depends on $p$ in an unexpectedly  regular way. It can be approximated up to a high power of $p$ 
 by a universal polynomial with large rational coefficients whose origin is not yet understood. Moreover, these coefficients again have a nice formula for the denominators.

\subsection{A \texorpdfstring{$q$}{q}-analog}

In the matrix (\ref{mdt}) for $D_t$ replace polynomial expressions in $n$ by $q$-analogs:
$$\begin{array}{c|cc}\, & x^{n-1}&x^n\\ \hline x^{n-1}& {q^{n-1}-1\over q-1}{q^n-1\over q-1}(t+1) & \big({q^n-1\over q-1}\big)^2 t\\
x^n & \big({q^n-1\over q-1}\big)^2 & {q^n-1\over q-1}{q^{n+1}-1\over q-1}(t+1)\end{array}$$
and multiply the operator by $(q-1)^2$. We obtain
$$\begin{array}{c|cc}\, & x^{n-1}&x^n\\ \hline x^{n-1}& (q^{n-1}-1)(q^n-1)(t+1) & (q^n-1)^2 t\\
x^n & (q^n-1)^2 & (q^n-1)(q^{n+1}-1)(t+1)\end{array}$$

Let replace everywhere  $q^n\leadsto q^{n+\epsilon}$, denote $\xi:=q^\epsilon$. The resulting matrix is
$$\begin{array}{c|cc}\, & x^{n-1+\epsilon}&x^{n+\epsilon}\\ \hline x^{n-1+\epsilon}& (q^{n-1}\xi-1)(q^n\xi-1)(t+1) & (q^n\xi-1)^2 t\\
x^{n+\epsilon} & (q^n\xi-1)^2 & (q^n\xi-1)(q^{n+1}\xi-1)(t+1)\end{array}$$

In terms of generators of quantum torus $X:x^n\mapsto x^{n+1},\,\,Y:x^n\mapsto q^n x^n$
this is operator
$$(Y-1)(qY-1)+(Y-1)^2 X+t\big[(Y-1)(qY-1)+X^{-1}(Y-1)^2\big]    $$
After replacement $X\leadsto t^{1/2}X$ we get:
$$ (Y-1)(qY-1)+t^{1/2}\big[(Y-1)^2 X+X^{-1}(Y-1)^2\big]+t (Y-1)(qY-1)   $$

There are two series $\xi_1(t),\xi_2(t)$ (``monodromies") for which there are solutions,
the characteristic polynomial is
$$(\xi-\xi_1(t))(\xi-\xi_2(t))=\xi^2-2\xi+1+C(t)\xi      $$
where  $$ C(t)=C_1 t+C_2 t^2+\dots =0\cdot t+\frac{-q^2+2q-1}{q^2+2q+1} t^{2} +\frac{q^{6} - 3 q^{5} + q^{4} + 2 q^{3} + q^{2} - 3 q + 1}{q^{6} + 5 q^{5} + 11 q^{4} + 14 q^{3} + 11 q^{2} + 5 q + 1}t^3+\dots=$$
$$=-{(q-1)^2\over (q+1)^2}t^2+{(q - 1)^{2}  (q^{4} - q^{3} - 2 q^{2} - q + 1)\over  (q^{2} + q + 1)(q + 1)^{4} } t^3-$$
$$-{ (q - 1)^{2}  (q^{10} - q^{9} + 5 q^{8} + 22 q^{7} + 47 q^{6} + 54 q^{5} + 47 q^{4} + 22 q^{3} + 5 q^{2} - q + 1)\over   (q^{2} + 1)  (q^{2} + q + 1)^{2}(q + 1)^{6}}t^4+\dots $$
and all coefficients lie in $\Z\big[q,q^{-1},\big({1\over q^k-1}\big)_{k=1,2,\dots}\big]$.

The denominators of $C_k,\,2\leqslant k\leqslant 14$ are

$k=2: \Phi_2^{\bf 2}$

$k=3: \Phi_3\cdot \Phi_2^{\bf 4}$

$k=4: \Phi_4\cdot \Phi_3^2\cdot \Phi_2^{\bf 6}$

$k=5: \Phi_5\cdot \Phi_4^2\cdot \Phi_3^3\cdot \Phi_2^{\bf 8}$

%$k=6: \Phi_6\cdot \Phi_5^2\cdot \Phi_4^3\cdot \Phi_3^4\cdot \Phi_2^{\bf 10}$

%$k=7: \Phi_7\cdot \Phi_6^2\cdot \Phi_5^3\cdot \Phi_4^4\cdot \Phi_3^5\cdot \Phi_2^{\bf 12}$

$\dots$

$k=14: \Phi_{14}\cdot \Phi_{13}^2\cdot \Phi_{12}^3\cdot\Phi_{11}^4\cdot \Phi_{10}^5\cdot\Phi_9^6\cdot\Phi_8^7\cdot\Phi_7^8\cdot\Phi_6^9\cdot\Phi_5^{10}\cdot\Phi_4^{11}\cdot\Phi_3^{12}\cdot \Phi_2^{\mathbf {26}}$
\newline
where $\Phi_n=q^{\varphi(n)}+\dots\in\Z[q]$ is the $n$-th cyclotomic polynomial. Presumably, this regular behavior $\prod_{l=3}^k \Phi_l^{k+1-l}\cdot \Phi_2^{2k-2}$ of the denominator can be proven using Theorems 5.2 and 5.4.

Expression ${C\over (q-1)^2}$   has no poles at $q=1$,  and its value at $q=1$ is exactly the series in the arithmetic case:
$$-{1\over 4} t^2-{1\over 24} t^3-{101\over 576 } t^4+\dots =\sum c_n t^n$$
Therefore, substituting $q\leadsto 1$ in the denominators of $C_k$ we get integers which are multiples of  the denominators of coefficients $c_k$. It is known that 
$$\Phi_n(q)=\prod_{d|n}(q^d-1)^{\mu(\frac{n}{d})}$$
which means that in the limit $q\to 1$ we can write 
$$\forall n\geqslant 2:\quad \Phi_n(1)=\prod_{d|n}d^{\mu(\frac{n}{d})}=\begin{cases}p&\text{ if }n=p^r\text{ is a power of prime}\\1&\text{ otherwise}\end{cases}.$$
Substituting these values in our hypothetical formula for denominators of $C_n$ we obtain a good approximation to  the hypothetical formulas (\ref{cn}) for denominators of $c_n$ in the arithmetic case. At least $p$-adic valuations in the case $\sqrt{n}<p\leqslant n$ match.

{\bf Remark 7.2.} For a given $n\geqslant 1$ let $\Q(\xi_n)\simeq\Q[q]/\Phi_n(q)$ be
the $n$-th cyclotomic field generated by a primitive root of unity $\xi_n$. 
The coefficients $C_m$ can be interpreted as elements of a non-archimedean field
$\Q(\xi_n)(\!(q-\xi_n)\!)$. Then the above formula for the denominators of $C_m$ 
indicates that the series $\sum_{m\geqslant 1}C_mt^m$ has a {\it positive} radius of convergence in the above field, in contrast to the arithmetic case.

 \addcontentsline{toc}{section}{Acknowledgements}

\section*{Acknowledgements}

  A.O. is grateful to IHES, where this work was done, for invitations and an excellent working atmosphere.

\end{document}